# Reconfiguring flexibility in renewable power-to-ammonia systems using molten-salt thermal energy storage in the ammonia synthesis loop: A coordinated electro-hydrogen-thermal scheduling approach

Yiwei Qiu[1], Qingjie Sun[1], Yangjun Zeng[1], Ge Chen[3], Longjie Yang[1], Ge He[2], Xu Ji[2], Shi Chen[1*], Buxiang Zhou[1], Kaigui Xie[1]

[1] College of Electrical Engineering, Sichuan University, Chengdu, 610065, China;
[2] School of Chemical Engineering, Sichuan University, Chengdu, 610065, China;
[3] School of Advanced Engineering, Great Bay University, Dongguan, 523808, China
*chen_shi@scu.edu.cn

**Abstract:** In renewable power-to-ammonia (ReP2A) systems, the intermittency of wind and solar generation propagates through electrolytic hydrogen production and induces thermal instability in the ammonia synthesis reactor (ASR). The resulting temperature cycling accelerates fatigue and shortens service life, while reactor thermal inertia limits flexible start-up, shutdown, and load adjustment. To address this issue, this study integrates molten-salt thermal energy storage (MS-TES) into the Haber-Bosch synthesis loop and develops a coordinated electro-hydrogen-thermal scheduling framework. MS-TES decouples hydrogen supply fluctuations from reactor thermal dynamics by enabling hot standby operation and sustained thermal support during start-up and low-load conditions. A state-space model is established to capture the thermal dynamics of the ASR and MS-TES. Based on this model, an optimal scheduling program coordinates ammonia synthesis operation with hydrogen production, battery energy storage (BES), and hydrogen storage (HS). The problem is formulated as a mixed-integer linear program (MILP) and extended with information gap decision theory (IGDT) to address renewable uncertainty. Case studies based on an industrial-scale project in northern China show that MS-TES enhances reactor thermal stability and system-level flexibility, while diminishing the marginal benefit of large BES capacity. As a result, a configuration combining small BES, HS, and MS-TES achieves near-equivalent performance to large-BES systems, with lower investment and improved economic returns. Year-round simulations further show that MS-TES avoids ASR start-up and shutdown and delivers consistently higher net revenue under variable renewable conditions.

| Nomenclature | Meaning |
|---|---|
| **Abbreviations** | |
| AS | Ammonia Synthesis |
| ASR | Ammonia Synthesis Reactor |
| BES | Battery Energy Storage |
| HB | Haber-Bosch |
| HS | Hydrogen Storage |
| IGDT | Information Gap Decision Theory |
| MILP | Mixed-Integer Linear Programming |
| MS | Molten Salt |
| MS-TES | Molten-Salt Thermal Energy Storage |
| PV | Photovaltaic |
| ReP2A | Renewable Power to Ammonia |
| RES | Renewable Energy Sources |
| RIG | Reactor Inlet Gas |
| ROG | Reactor Outlet Gas |
| SUH | Start-Up Heater |
| WT | Wind Turbine |
| | |
| **Indicators** | |
| $t$ | Indicator of time intervals |
| $\Delta t$ | Indicator of step length |
| $\tau$ | Scheduling horizon |
| | |
| **Variable** | |
| $b_t^{AS,on/by/cs/off}$ | AS operating state of production/hot standby/cold start/shutdown |
| $b_t^{AS,startup/shutdown}$ | Action of the AS section entering/exiting the production state |
| $b_t^{AS,Inoff/Outoff}$ | Action of the AS section entering/exiting the shutdown state |
| $b_t^{MS,on/by}$ | MS-TES gas-heating/standby state |
| $b_t^{BAT,cha/dis}$ | BES charging/discharging state |
| $E^{BAT}$ | BES stored energy |
| $F^{ASR,in/out}$ | RIG/ROG flow rate |
| $F^{ASR,in,direct}$ | Gas flow rate entering the ASR directly |
| $F^{ASR,H_2}$ | Hydrogen consumption rate |
| $F^{ASR,NH_3}$ | Ammonia production rate |
| $F^{HP,H_2}$ | Hydrogen production rate |
| $F^{H_2,out}$ | Hydrogen outflow of hydrogen storage tank |
| $F^{MS}$ | Gas flow rate heated by MS-TES |
| $F^{SU}$ | Gas flow rate heated by the SUH |
| $L$ | AS load level |
| $P^{ASR,aux}$ | Power consumption of AS |
| $P^{BAT,cha/dis}$ | BES charging/discharging power |
| $P^{HP}$ | Power for hydrogen production |
| $P^{Grid}$ | Power purchase from the grid |
| $Q^{ASR,react}$ | ASR reaction heat duty |
| $Q^{ASR,in/out}$ | RIG/ROG heat duty |
| $Q^{ASR,cool}$ | ASR cooling duty provided by quenching or heat exchange |
| $Q^{ASR,diss}$ | ASR heat duty dissipated to the ambient |
| $Q^{MS,heat}$ | MS-TES heat duty supplied to the RIG |
| $Q^{MS,diss}$ | Heat duty dissipated from the MS-TES |
| $Q^{SU,heat}$ | Heat duty supplied by the SUH to the RIG |
| $T^{ASR}$ | ASR temperature |
| $T^{ASR,in/out}$ | RIG/ROG temperature |
| $T^{MS,ms}$ | MS temperature |
| $T^{MS,out}$ | RIG temperature after MS-TES heating |
| $T^{SU,out}$ | RIG temperature after SUH heating |
| $\alpha$ | Uncertainty horizon of RES |
| | |
| **Parameters** | |
| $c^{ASR,in/out}$ | RIG/ROG specific heat capacity |
| $c^{ms}$ | MS specific heat capacity |
| $c^{HP}$ | Power consumption per unit hydrogen production |
| $c^{ASR,aux}$ | AS auxiliary power consumption per unit load rate |
| $C^{ASR}$ | ASR lumped thermal capacitance |
| $F^{ASR,in,by}$ | RIG flow rate in hot standby state |
| $F^{ASR,in,0}$ | RIG flow rate intercept |
| $F^{ASR,in,1}$ | RIG flow rate slope coefficient |
| $\overline{F}^{HP,H_2}$ | Limit of the hydrogen production rate |
| $\overline{F}^{MS}$ | Maximum gas flow rate through the MS-TES heater |
| $\overline{F}^{SU}$ | Maximum gas flow rate through the SUH |
| $T^{am}$ | Ambient temperature |
| $R^{ASR,diss}$ | ASR equivalent thermal resistance to ambient |
| $R^{MS,diss}$ | MS-TES equivalent thermal resistance to ambient |
| $\underline{T}^{ASR,act}, \overline{T}^{ASR,act}$ | ASR temperature limits in production state |
| $\underline{T}^{MS,ms}, \overline{T}^{MS,ms}$ | MS temperature limits |
| $\overline{T}^{MS,out}$ | Upper limit of the gas temperature after MS-TES heating |
| $\overline{T}^{SU,out}$ | Upper limit of the gas temperature after SUH heating |
| $P^{ASR,aux,0}$ | Baseline power of AS auxiliary |
| $\overline{P}^{MS,heat}$ | Maximum electric heating power of the MS-TES electric heater |
| $\overline{P}^{SU,heat}$ | Maximum heating power of the SUH |
| $\overline{P}^{BAT,cha/dis}$ | BES charging/discharging power limit |
| $\underline{n}^{H_2,Buffer}, \overline{n}^{H_2,Buffer}$ | Hydrogen inventory limits |
| $\underline{L}, \overline{L}$ | AS operational load limits |
| $\underline{r}, \overline{r}$ | AS ramping limits |
| $\overline{Q}^{ASR,cool}$ | Upper limit of the ASR cooling duty |
| $\overline{Q}^{MS,heat}$ | Upper limit of the MS-TES heating duty |
| $\overline{Q}^{SU,heat}$ | Upper limit of the SUH heating duty |
| $S^{ASR,H_2}$ | Hydrogen consumption at rated load |
| $S^{ASR,NH_3}$ | Ammonia production at rated load |

# 1. Introduction

## 1.1. Background and motivation

As the global energy transition accelerates, the renewable power-hydrogen-ammonia pathway is increasingly recognized as a promising option for improving renewable energy utilization and large-scale, long-duration energy storage [1, 2]. Green hydrogen production provides flexible electrical demand that can respond rapidly to renewable variability [3-6] and supports decarbonization across multiple sectors [7-9]. Ammonia synthesis (AS), traditionally a major hydrogen consumer, is also emerging as an energy carrier [2, 10]. In China alone, more than 100 green hydrogen and ammonia projects are planned, corresponding to an estimated annual green ammonia capacity of 15 million tonnes [11]. Similar developments are underway in Europe, South Asia, the Middle East, etc. [12].

Conventional Haber-Bosch AS processes are designed for steady high-temperature and high-pressure operation and are therefore sensitive to load changes [13]. In renewable power-to-ammonia (ReP2A) systems, fluctuations in wind and solar generation propagate through hydrogen production and induce load variations in the AS unit, leading to temperature fluctuations in the ammonia synthesis reactor (ASR). Under prolonged low-renewable conditions, the AS unit may shut down. Repeated thermal cycling accelerates catalyst degradation and structural fatigue, reducing reliability and service life [14-16], while reactor thermal inertia limits flexible start-up, shutdown, and load adjustment. Maintaining ASR temperature stability under intermittent hydrogen supply thus remains a key challenge.

Existing studies have explored several approaches to improve adaptability to renewable intermittency, which include: 1) deploying battery energy storage (BES) to mitigate short-term power fluctuations [17]; 2) using hydrogen storage (HS) for intertemporal energy shifting [17-19]; and 3) enhancing the operational flexibility of the AS process to enable dynamic load regulation [15, 20-27].

However, these approaches primarily address electrical or material buffering and do not directly resolve reactor thermal stress. As a result, renewable intermittency can still induce significant ASR temperature fluctuations. Expanding BES and HS can mitigate these effects, but it leads to high capital costs and may be constrained by safety regulations [17, 19].

A complementary strategy is to enable hot standby operation, which maintains the ASR within an appropriate temperature range during hydrogen shortages. This approach directly addresses the thermal behavior of the reactor, which cannot be effectively managed through electrical or hydrogen storage alone, reduces start-up and shutdown cycles, and limits thermomechanical stress.

In this context, molten-salt thermal energy storage (MS-TES) has emerged as a promising solution. Demonstration projects in Otog Banner, Inner Mongolia [28], and Songyuan, Jilin [29], in northern China, have integrated MS-TES into the AS loop, where surplus electricity is converted into stored heat and later released to preheat the reactor inlet gas (RIG) during hydrogen shortages. This decouples thermal management from instantaneous hydrogen supply and enables sustained reactor temperature control.

Although MS-TES is a mature technology widely applied in concentrating solar power (CSP) [30] and thermal power systems [31, 32], its integration into AS introduces coupled thermal-operational dynamics that are not captured in existing studies. In particular, the interaction between reactor thermal inertia, storage dynamics, and system-level scheduling decisions remains insufficiently understood.

To address these gaps, this study develops a dynamic thermal model of the MS-TES-integrated AS process and an optimization-based scheduling framework for ReP2A systems. With coordinated electro-hydrogen-thermal flexibility, the proposed approach provides a basis for safe, efficient, and cost-effective operation. Section 1.2 reviews related work, and Section 1.3 summarizes the main contributions.

## 1.2. Literature review

To accommodate renewable variability and ensure stable ASR operation, existing research on ReP2A systems can be broadly grouped into three categories: deployment of BES and HS facilities, optimization-based scheduling and control, and process-level improvements to AS.

**a) Deployment of BES and HS.** BES and HS differ in response time, energy density, cost, and cycle life. BES provides fast regulation over timescales from seconds to hours, while HS enables long-duration energy shifting [17, 18]. Hybrid storage systems have been proposed to address flexibility requirements across multiple timescales [33-34].

Table 1 Summary of studies on ReP2A system modeling and scheduling

| Literature | Components | Flexibility of the AS | | | Heat management | | Renewable uncertainty | Method |
|---|---|---|---|---|---|---|---|---|
| | | Load range limits | Ramping limits | Standby ability | Reactor heat dynamics | Thermal storage | | |
| Wu et al. [22] | PV[1]/WT[2]/BES/EL[3]/HS/AS | √ | √ | × | × | × | × | MILP |
| Wang [23] | PV/WT/BES/EL/HS/AS | √ | √ | × | × | × | × | Two-stage stochastic MILP |
| Shi et al. [24] | PV/WT/EL/HS/AS | √ | √ | × | × | × | √ | Chance-constrained programming |
| Zhou et al. [33] | PV/WT/GE[4]/EL/HS/AS | √ | √ | × | × | × | √ | MILP |
| Yu et al. [37] | PV/WT/EL/HS/AS | √ | √ | × | × | × | √ | Two-stage optimization |
| Zhou et al. [38] | PV/WT/EL/HS/AS | √ | √ | × | × | × | × | PSO-MILP |
| Wu et al. [39] | PV/WT/EL/HS/AS | √ | √ | × | × | × | √ | Two-stage robust optimization |
| Yu et al. [40] | PV/WT/BES/EL/HS/AS | √ | √ | × | × | × | √ | IGDT-MILFP |
| Allman et al. [41] | WT/EL/HS/AS | √ | √ | × | × | × | × | Receding horizon optimization |
| Li et al. [42] | PV/WT/BES/EL/HS/AS | √ | √ | × | × | × | × | MIQP |
| Wu et al. [43] | PV/WT/EL/HS/AS | √ | √ | √ | × | × | × | DR[5]-guided control |
| Ding et al. [44], Zheng et al. [45] | PV/WT/EL/HS/AS | √ | √ | × | × | × | √ | MILP |
| **This work** | **PV/WT/BES/EL/HS/AS/MS-TES** | √ | √ | √ | √ | √ | √ | IGDT-MILP |

Note:[1] Photovoltaic; [2] Wind turbine; [3] Electrolyzer; [4] Gravity energy storage; [5] Dispatchable region; √ Considered; × Not considered.

In principle, sufficiently large BES and HS can smooth renewable variability and stabilize hydrogen supply. In practice, deployment is constrained by cost and safety. BES requires high capital investment [17], and HS introduces additional costs related to compression, safety compliance, and efficiency losses [19]. Safety risks also increase significantly when HS capacity exceeds approximately $10^5$ Nm$^3$ [35]. According to the Chinese standard GB/T 29729-2022, *Essential requirements for the safety of hydrogen systems*, stricter siting and emergency control requirements apply to large-scale hydrogen systems. For an AS plant producing $2\times10^5$ t/yr of ammonia, an HS capacity of $5\times10^5$ Nm$^3$ supports full-load operation for only several to just over ten hours. Consequently, BES and HS alone are insufficient to bridge multi-day periods of low wind and solar output.

**b) Scheduling and control for enhanced flexibility.** Studies have focused on improving system-level flexibility through improved scheduling and control. Rosbo et al. [20] combined dynamic AS modeling with stability analysis to expand the feasible load range. Kong et al. [21] applied nonlinear model predictive control to improve temperature and pressure tracking in Haber-Bosch processes. Wu et al. [22] developed a multi-timescale capacity planning and scheduling model to maximize annual profit, while Wang [23] proposed a stochastic two-stage MILP model solved using improved Benders decomposition. Shi et al. [24] employed chance-constrained scheduling for coordinated operation of power-hydrogen-ammonia networks.

Although these methods improve operational flexibility, they typically treat the AS process as a load-following unit and do not explicitly model reactor thermal dynamics or its interaction with scheduling decisions. As a result, hydrogen supply variability continues to propagate into ASR temperature fluctuations, accelerating catalyst degradation [16] and structural fatigue [36].

**c) Process-level improvements to AS.** Other studies have sought to enhance AS flexibility through reactor and process design. Fahr et al. [15] developed a detailed dynamic model of the Haber-Bosch process and optimized ASR design and pressure control. Ji et al. [25] proposed a multi-steady-state process to better accommodate renewable variability. Additional thermodynamic and design optimizations have expanded the load range [26, 27]. However, these approaches focus on process design and local control, and do not address system-level variability in hydrogen supply, making sustained temperature stability difficult under highly intermittent conditions.

Integrating MS-TES into the AS loop provides a direct mechanism to regulate reactor thermal inertia. During hydrogen shortages, stored thermal energy maintains hot standby operation, keeping the ASR within the catalyst's active range. Once the hydrogen supply recovers, the system can quickly resume production, reducing start-up delays and thermal stress.

Despite its maturity in CSP and conventional power systems [30-32], MS-TES has received limited attention in ReP2A systems. Existing studies rarely capture the coupled dynamics between thermal storage, reactor thermal inertia, and system-level scheduling. As a result, the role of thermal storage in coordinating multi-energy flexibility remains unclear. This study addresses this gap by developing a dynamic model and scheduling framework for MS-TES-integrated ReP2A systems. Related studies are summarized in Table 1.

### *1.3. Contributions of this work*

To improve the flexibility and economic performance of ReP2A systems under intermittent power supply, this study develops an MS-TES-integrated scheduling framework. The main contributions are as follows:

1) A state-space thermal model is developed to capture the coupled heat dynamics of the ASR and MS-TES. The model preserves key thermal inertia and heat-loss characteristics while remaining compatible with system-level optimization, enabling reactor temperature constraints to be explicitly embedded into a tractable scheduling problem.

2) The AS thermal dynamics are embedded into an MILP-based scheduling model that co-optimizes AS operating states (*shutdown*, *cold start*, *hot standby*, and *production*) with hydrogen production, BES, HS, and renewable supply. Renewable uncertainty is addressed using information gap decision theory (IGDT).

3) Simulations based on a real-life project in northern China show that MS-TES improves ASR temperature stability. In particular, long-duration thermal support reduces reliance on large BES capacity while maintaining stable ASR operation. A coordinated configuration combining small BES, HS, and MS-TES achieves near-equivalent production performance to large-BES systems, while reducing storage investment and improving net revenue.

The remainder is organized as follows. Section 2 presents the dynamic thermal model of the MS-TES-integrated AS process. Section 3 formulates the optimal scheduling model for the ReP2A system. Section 4 reports and discusses the case study results. Section 5 concludes the paper.

## 2. Dynamic thermal modeling of MS-TES-integrated AS process
### *2.1. Configuration of the utility-scale ReP2A system*

The ReP2A system comprises wind and PV generation, water electrolysis for hydrogen production, air separation, compression and buffering, and AS [46]. The AS section adopts the Haber-Bosch process and includes the ASR, ammonia separation, and a recycle gas loop [12]. In this study, we integrate an MS-TES module in the loop. The overall configuration is shown in Fig. 1.

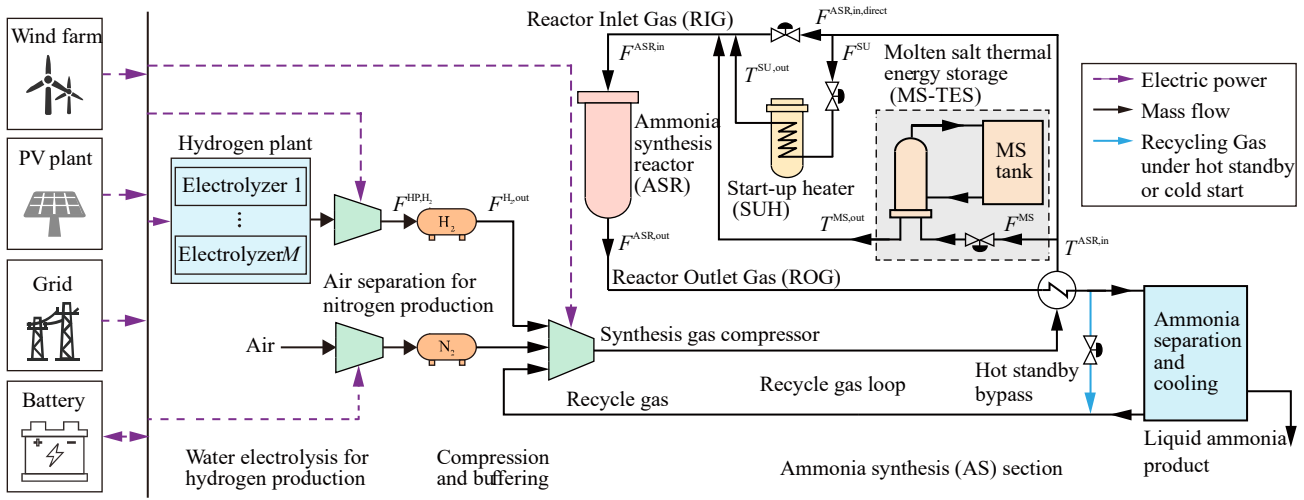

**Fig. 1**. Schematic diagram of the ReP2A system with an MS-TES-integrated AS process.

## 2.2. Overview of the AS process integrating MS-TES

As shown in Fig. 1, a single-tank MS-TES module is integrated to the RIG side of the AS process. During low renewable output and insufficient hydrogen supply, the AS can switch to hot standby. The recycle gas bypasses the ammonia cooling and separation and is directly routed back to the compressor, indicated by the blue path in Fig. 1. The gas is preheated by the start-up heater (SUH) and/or MS-TES before entering the ASR. This compensates for heat losses and maintains the catalyst within its active temperature range. This enables rapid restart when hydrogen supply recovers and reduces thermal stress. The same heating approach is applied under low-load conditions when reaction heat is insufficient to sustain reactor temperature.

After releasing heat, the molten salt (MS) must be reheated electrically to restore its temperature and avoid solidification. Heat can be supplied to the RIG only when the MS temperature exceeds the gas temperature. If hot standby is required when the MS temperature is low, the SUH provides direct electrical heating. Although effective, SUH does not store thermal energy. The scheduling program therefore needs to coordinate MS-TES and SUH.

AS can adjust its load to follow hydrogen supply fluctuations. The load level, denoted by $L$, typically varies from 30% to 110%, with a ramping limit of ±25%/h [21]. The AS operates in four states, i.e., *production*, *shutdown*, *cold start*, and *hot standby*, shown in Fig. 2. State transitions are detailed in Section 3.1.

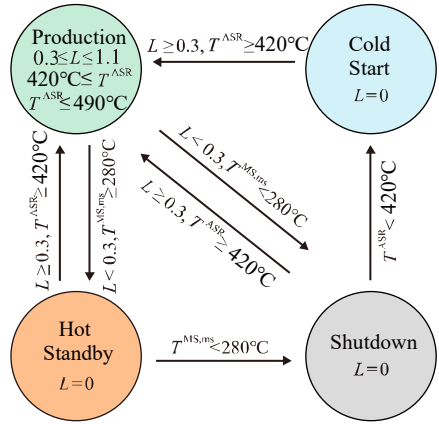

**Fig. 2**. Operating states and transitions of the AS process with hot-standby ability.

Fig. 3 illustrates the evolution of ASR load, temperature, and RIG heating duty across operating states. During shutdown from $t_1$ to $t_2$, the ASR temperature decays due to heat loss. During cold start from $t_2$ to $t_3$, the SUH and MS-TES heat the RIG to raise the temperature until it reaches the lower limit for production. In hot standby from $t_4$ to $t_5$, RIG heating maintains ASR temperature within the allowable range.

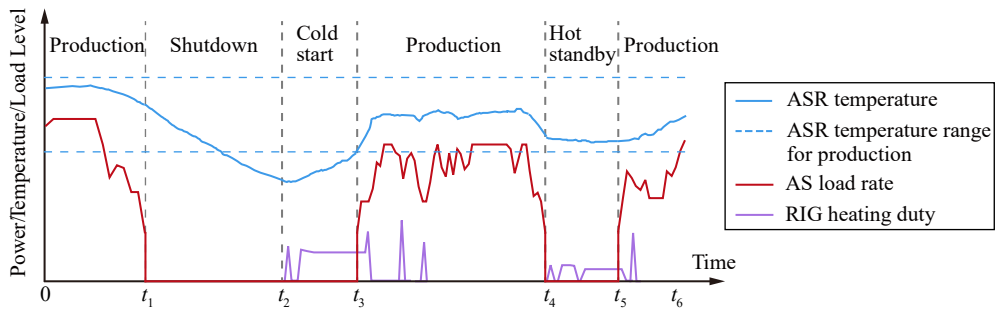

**Fig. 3**. Illustration of load and temperature dynamics under state transitions of the AS process integrating MS-TES.

## 2.3. Dynamic model of the ASR temperature

As this study focuses on system-level scheduling, the ASR is represented by a lumped-parameter thermal model, which follows

$$C^{ASR} \frac{d}{dt} T^{ASR} = Q^{ASR,react} + Q^{ASR,in} - Q^{ASR,out} - Q^{ASR,cool} - Q^{ASR,diss}, \quad (1)$$

where $T^{ASR}$ is the ASR temperature; $C^{ASR}$ is the lumped thermal capacitance, including the reactor structure, the catalyst bed, and internal gas stream; $Q^{ASR,react}$ is the reaction heat; $Q^{ASR,in}$ and $Q^{ASR,out}$ are the heat duties of RIG and reactor outlet gas (ROG), respectively; $Q^{\text{ASR,cool}}$ is the cooling duty; $Q^{\text{ASR,diss}}$ denotes the heat dissipated to the ambient, which dominates ASR temperature decay after shutdown.

Terms in (1) are modeled as follows.

**a) Heat loss.** Heat transfer to the ambient is considered, as

$$Q^{\text{ASR,diss}} = \frac{T^{\text{ASR}} - T^{\text{am}}}{R^{\text{ASR,diss}}}, \tag{2}$$

where $T^{\text{am}}$ is the ambient temperature; $R^{\text{ASR,diss}}$ is the thermal resistance between the ASR and the environment.

**b) Reaction heat.** Reaction heat is generated only in the production state and depends on the load rate $L$, as $Q^{\text{ASR,react}} = Q^{\text{ASR,react}}(L)$. When the AS is not producing, $Q^{\text{ASR,react}} = 0$.

The AS reaction typically occurs at 400-500 °C and 150-200 bar in the presence of an iron-based catalyst [40]:

$$\frac{3}{2}\text{H}_2(g) + \frac{1}{2}\text{N}_2(g) \to \text{NH}_3(g), \Delta H \approx -52.3 \text{ kJ/mol}. \tag{3}$$

The reaction heat is approximated as a linear function of load $L$:

$$Q^{\text{ASR,react}}(L) = c^{Q,\text{react}} L, \tag{4}$$

where $c^{Q,\text{react}}$ is a constant factor. Although nonlinear reaction kinetics exothermic models [25] could be adopted and handled by piecewise approximations, the resulting improvement in accuracy is limited and therefore omitted.

**c) Reactor inlet and outlet heat duties.** In the production state, the heat duties of the RIG and ROG depend on their flow rates, temperatures, and specific heat capacities:

$$Q^{\text{ASR,in/out}} = F^{\text{ASR,in/out}} c^{\text{ASR,in/out}} T^{\text{ASR,in/out}}, \tag{5}$$

where $F^{\text{ASR,in/out}}$, $c^{\text{ASR,in/out}}$, and $T^{\text{ASR,in/out}}$ denote the gas flow, specific heat capacity, and temperature, respectively.

In practice, the AS design package specifies RIG and ROG temperatures based on plant scale, flowsheets, and catalyst. During operation, the control system maintains these variables near their design values across the load range. Design data for a 200 kt/yr AS plant, shown in Fig. A1 in Appendix A, indicate that RIG and ROG temperatures and specific heat capacities remain approximately constant between 50% and 110% load, while flow rates vary linearly with $L$. Accordingly, in the production state, we have

$$T^{\text{ASR,in/out}} = \tilde{T}^{\text{ASR,in/out}}, \ c^{\text{ASR,in/out}} = \tilde{c}^{\text{ASR,in/out}}, \ F^{\text{ASR,in/out}} = \tilde{F}^{\text{ASR,in/out}}(L), \tag{6}$$

where $\tilde{T}^{\text{ASR,in/out}}(L)$, $\tilde{c}^{\text{ASR,in/out}}(L)$, and $\tilde{F}^{\text{ASR,in/out}}(L)$ are the process design values of the RIG/ROG temperature, specific heat capacity, and flow rate, respectively.

In contrast, in hot standby, the ROG bypasses the separation unit and is directly routed back to the compressor (blue path in Fig. 1). The recycle gas temperature is actively controlled to a constant value higher than that in production, denoted by $\tilde{T}^{\text{ASR,in,by}}$, say $T^{\text{ASR,in}} = \tilde{T}^{\text{ASR,in,by}}$, to compensate for ASR heat dissipation.

**d) Cooling duty.** By regulating quench gas flow (for quench-type reactors) or heat-exchange gas flow (for heat-exchange-type reactors), ASR temperature can be actively controlled [47]. The cooling duty $Q^{\text{ASR,cool}}$ is therefore treated as a control variable.

### 2.4. Thermal management during hot standby and cold start

Fig. 2 illustrates the transition logic among four states. When reaction heat is sufficient, the RIG enters the ASR without preheating. During hot standby or cold start, reaction heat is absent, and the RIG must be heated before entering the reactor to raise its temperature.

In these states, the RIG is divided into three streams, with one entering the ASR directly, one heated by the MS-TES, and one heated by the SUH, as

$$F^{\text{ASR,in}} = F^{\text{ASR,in,direct}} + F^{\text{MS}} + F^{\text{SU}}, \tag{7}$$

where $F^{\text{ASR,in,direct}}$, $F^{\text{MS}}$, and $F^{\text{SU}}$ are the corresponding gas streams. The resulting RIG heat duty is:

$$Q^{\text{ASR,in}} = c^{\text{ASR,in}} \left( F^{\text{ASR,in,direct}} T^{\text{ASR,in}} + F^{\text{MS}} T^{\text{MS,out}} + F^{\text{SU}} T^{\text{SU,out}} \right), \tag{8}$$

where $T^{\text{MS,out}}$ and $T^{\text{SU,out}}$ are the gas temperatures after being heated by MS-TES and SUH, respectively.

**a) Molten salt thermal energy storage (MS-TES).** A binary molten salt composed of 40% KNO$_3$ and 60% NaNO$_3$, commonly used in CSP plants, is adopted in the AS process. It has a melting point of approximately 220 °C, remains chemically stable up to 600 °C, and typically operates between 280 °C and 565 °C [48].

When activated, MS is pumped through a heat exchanger to transfer heat to the RIG. The heated gas flow $F^{MS}$ is controlled by a valve. Because MS temperature varies much more slowly than the gas residence time in the exchanger, we assume the outlet gas temperature $T^{MS,out}$ is equal to the MS temperature minus a small fixed gap $\Delta T^{MS}$, i.e., $T^{MS,out} = T^{MS,ms} - \Delta T^{MS}$.

The MS temperature is governed by:

$$r^{ms}V^{ms}c^{ms}\frac{d}{dt}T^{MS,ms} = -\frac{1}{\eta^{MS,gas,heat}}Q^{MS,heat} + \eta^{MS,heat}P^{MS,heat} - Q^{MS,diss}, \tag{9}$$

where $r^{ms}$, $V^{ms}$, and $c^{ms}$ are the density, volume, and specific heat capacity of the MS; $P^{MS,heat}$ is the electric heating power; $\eta^{MS,gas,heat}$ and $\eta^{MS,heat}$ are the efficiencies of the heat exchanger and electric heater; $Q^{MS,heat}$ is the heat transferred to the RIG, governed by

$$Q^{MS,heat} = c^{ASR,in}F^{MS}\left(T^{MS,out} - T^{ASR,in}\right), \tag{10}$$

and $Q^{MS,diss}$ is the heat dissipated to the ambient:

$$Q^{MS,diss} = \frac{T^{MS,ms} - T^{am}}{R^{MS,diss}}, \tag{11}$$

where $R^{MS,diss}$ denotes the thermal resistance between MS-TES and the environment.

**b) Start-up heater (SUH).** The SUH provides direct electric heating to the RIG:

$$\eta^{SU,heat}P^{SU,heat} = Q^{SU,heat} = c^{ASR,in}F^{SU}\left(T^{SU,out} - T^{ASR,in}\right), \tag{12}$$

where $P^{SU,heat}$ and $\eta^{SU,heat}$ are the electric power and efficiency of SUH; $Q^{SU,heat}$ is the heat duty delivered to the RIG.

## 3. Optimal scheduling of the ReP2A system integrating MS-TES

Based on the models developed in Section 2, we formulate an optimal scheduling program in which reactor temperature evolution is explicitly coupled with operating decisions. The objective is to maximize net revenue while ensuring safe operation. At each time step, the model jointly determines hydrogen production, AS operating state and load, and the operation of MS-TES and SUH. For convenience, we use subscript $t$ to denote variables at time step $t$, $\Delta t$ the step length, and $\tau$ the scheduling horizon.

### 3.1. AS operating state transition logic

As shown in Fig. 2, AS operates in four mutually exclusive states, i.e., *production*, *hot standby*, *cold start*, and *shutdown*. This is enforced by

$$b_t^{AS,on} + b_t^{AS,by} + b_t^{AS,cs} + b_t^{AS,off} = 1, \tag{13}$$

$$0 \leq b_t^{AS,on}, b_t^{AS,by}, b_t^{AS,cs}, b_t^{AS,off} \leq 1, \tag{14}$$

where $b_t^{AS,on}, b_t^{AS,by}, b_t^{AS,cs}$, and $b_t^{AS,off}$ are binary variables indicating the four states.

Binary variables $b_t^{AS,startup}$ and $b_t^{AS,shutdown}$ denote transitions into and out of production. The logic is enforced by

$$b_t^{AS,on} - b_{t-1}^{AS,on} = b_t^{AS,startup} - b_t^{AS,shutdown}, \tag{15}$$

$$b_t^{AS,startup} + b_t^{AS,shutdown} \leq 1. \tag{16}$$

Similarly, entry into and out of shutdown are represented by $b_t^{AS,Inoff}$ and $b_t^{AS,Outoff}$:

$$b_t^{AS,off} - b_{t-1}^{AS,off} = b_t^{AS,Inoff} - b_t^{AS,Outoff}, \tag{17}$$

$$b_t^{AS,Inoff} + b_t^{AS,Outoff} \leq 1. \tag{18}$$

As exemplified in Fig. 3, entry into production is allowed only when the ASR temperature activates the catalyst when $T^{ASR} \geq \underline{T}^{ASR,act}$, i.e.,

$$\begin{cases} \underline{T}^{\text{ASR,act}} \leq T_t^{\text{ASR}}, & b_t^{\text{AS,on}} = 1 \\ -\infty \leq T_t^{\text{ASR}}, & b_t^{\text{AS,on}} = 0 \end{cases}, \tag{19}$$

where $\underline{T}^{\text{ASR,act}}$ is the minimum allowable ASR temperature for production.

To incorporate this condition into a mixed integer linear program (MILP), we apply the Big-M method [49], as

$$\underline{T}^{\text{ASR,act}} - M\left(1 - b_t^{\text{AS,on}}\right) \leq T_t^{\text{ASR}}, \tag{20}$$

with $M$ being a sufficiently large constant.

In addition, the AS load level $L$ is bounded and subject to ramping constraints [25-27]:

$$\underline{L} b_t^{\text{AS,on}} \leq L_t \leq \overline{L} b_t^{\text{AS,on}}, \tag{21}$$

$$\underline{r} - M\left(b_t^{\text{AS,startup}} + b_t^{\text{AS,shutdown}}\right) \leq L_t - L_{t-1} \leq \overline{r} + M\left(b_t^{\text{AS,startup}} + b_t^{\text{AS,shutdown}}\right), \tag{22}$$

$$\underline{L} - M\left(1 - b_{t-1}^{\text{AS,startup}} - b_t^{\text{AS,shutdown}}\right) \leq L_{t-1} \leq \underline{L} + M\left(1 - b_{t-1}^{\text{AS,startup}} - b_t^{\text{AS,shutdown}}\right), \tag{23}$$

where $\overline{L}$ and $\underline{L}$ are the load level limits; $\overline{r}$ and $\underline{r}$ are ramping limits; (22) enforces smooth load change; and (23) fixes the load to $\underline{L}$ at start-up and immediately before shutdown.

### 3.2. Operational and mass-flow constraints of the AS process

**a) Molten salt thermal energy storage (MS-TES).** Similar to AS state transitions, binary $b_t^{\text{MS,on}}$ indicates heating state, and $b_t^{\text{MS,by}}$ indicates standby when the MS-TES module does not supply heat to the RIG but the MS can be electrically heated. The logic follows

$$b_t^{\text{MS,on}} + b_t^{\text{MS,by}} = 1, \tag{24}$$

$$b_t^{\text{MS,on}} - \left(b_t^{\text{AS,on}} + b_t^{\text{AS,by}} + b_t^{\text{AS,cs}}\right) \leq 0, \tag{25}$$

$$0 \leq b_t^{\text{MS,on}}, b_t^{\text{MS,by}} \leq 1, \tag{26}$$

where (25) means that the MS-TES can provide heat only when the AS is in production, hot standby, or cold start. The heated gas flow is limited by

$$0 \leq F_t^{\text{MS}} \leq b_t^{\text{MS,on}} \overline{F}^{\text{MS}}, \tag{27}$$

where $\overline{F}^{\text{MS}}$ is the maximum allowable flow.

MS temperature is constrained to avoid freezing and overheating:

$$\underline{T}^{\text{MS,ms}} \leq T_t^{\text{MS,ms}} \leq \overline{T}^{\text{MS,ms}}, \tag{28}$$

$$T^{\text{ASR,in}} - M\left(1 - b_t^{\text{MS,on}}\right) \leq T_t^{\text{MS,ms}}, \tag{29}$$

where the latter ensures that heat transfer occurs only when MS is hotter than the inlet gas. The Big-M formulation resembles (21).

The electric heater's power is bounded by

$$0 \leq P_t^{\text{MS,heat}} \leq \overline{P}^{\text{MS,heat}}. \tag{30}$$

**b) Start-up heater (SUH).** The SUH operates when the AS is in hot standby, cold start, and low load conditions, and is inactive in shutdown. Gas flow, electrical power, and heater outlet gas temperature satisfy:

$$0 \leq F_t^{\text{SU}} \leq \left(1 - b_t^{\text{AS,off}}\right) \overline{F}^{\text{SU}}, \tag{31}$$

$$0 \leq P_t^{\text{SU,heat}} \leq \left(1 - b_t^{\text{AS,off}}\right) \overline{P}^{\text{SU,heat}}, \tag{32}$$

$$T^{\text{ASR,in}} \leq T_t^{\text{SU,out}} \leq \overline{T}^{\text{SU,out}}, \tag{33}$$

where $\overline{F}^{\text{SU}}$, $\overline{P}^{\text{SU,heat}}$, and $\overline{T}^{\text{SU,out}}$ are the limits of the gas flow, heating power, and outlet temperature.

**c) ASR cooling.** Excess reaction heat is removed through quenching or heat exchange. Cooling duty is automatically regulated to stabilize temperature [25], following

$$0 \leq Q_t^{\text{ASR,cool}} \leq b_t^{\text{AS,on}} \overline{Q}^{\text{ASR,cool}}, \tag{34}$$

where $\overline{Q}^{\text{ASR,cool}}$ is the maximum cooling capacity.

**d) Reactor gas flow and temperature control.** Since MS-TES and SUH draw gas from the RIG, as shown in Fig. 1, we have:

$$F_t^{MS} + F_t^{SU} \leq F^{ASR,in}. \tag{35}$$

In summary, the total RIG flow rate $F^{ASR,in}$ depends on the AS operation state. During production ($b_t^{AS,on} = 1$), the $F^{ASR,in}$ depends on the load level $L$, as shown in Fig. A1. During hot standby and cold start ($b_t^{AS,by} = 1$ and $b_t^{AS,cs} = 1$), $F^{ASR,in}$ is set to a constant $F^{ASR,in,by}$ to deliver heat from MS-TES and SUH to the ASR. When shutdown ($b_t^{AS,off} = 1$), the RIG flow is set to zero. Accordingly, the total RIG flow rate satisfies:

$$F_t^{ASR,in} = b_t^{AS,on} F^{ASR,in}(L_t) + (b_t^{AS,by} + b_t^{AS,cs}) F^{ASR,in,by}. \tag{36}$$

Mass conservation requires

$$F_t^{ASR,out} = F_t^{ASR,in}. \tag{37}$$

Finally, hydrogen consumption $F^{ASR,H_2}$ and ammonia production $F^{ASR,NH_3}$ scale linearly with load $L$:

$$F_t^{ASR,H_2} = L_t S^{ASR,H_2}, \tag{38}$$

$$F_t^{ASR,NH_3} = L_t S^{ASR,NH_3}, \tag{39}$$

where $S^{ASR,H_2}$ and $S^{ASR,NH_3}$ are rated values corresponding to $L=1$.

Following discussions in Section 2.3c), the RIG and ROG temperatures are also associated with the AS operation state and given by

$$T^{ASR,in} = b_t^{AS,on} \tilde{T}^{ASR,in} + (b_t^{AS,by} + b_t^{AS,cs}) \tilde{T}^{ASR,in,by}, \tag{40}$$

$$T^{ASR,out} = b_t^{AS,on} \tilde{T}^{ASR,out} + (b_t^{AS,by} + b_t^{AS,cs}) \tilde{T}^{ASR,in,by}. \tag{41}$$

**e) Electric power demand of the AS.** Electricity consumed by compressors and other auxiliaries (excluding MS-TES and SUH), denoted by $P_t^{ASR,aux}$, is modeled as linear to the load level $L$ [41]:

$$P_t^{ASR,aux} = P^{ASR,aux,0}(b_t^{AS,on} + b_t^{AS,by} + b_t^{AS,cs}) + c^{ASR,aux} L_t, \tag{42}$$

where $P^{ASR,aux,0}$ and $c^{ASR,aux}$ are the baseline power and coefficient.

### 3.3. Reformulation of nonlinearities

The original heat-duty expressions in (5), (8), (10), and (12) introduce bilinear terms that prevent direct integration into an MILP framework. To preserve tractability, we reformulate these terms using energy-based representations, which replace temperature-flow products with bounded heat variables, enabling the thermal dynamics to be embedded into a linear scheduling model while retaining physical consistency.

Specifically, the RIG heat duty $Q_t^{ASR,in}$ is expressed directly as

$$\begin{aligned} Q_t^{ASR,in} &= c^{ASR,in}\left(F_t^{ASR,in,direct} T^{ASR,in} + F_t^{MS} T_t^{MS,out} + F_t^{SU} T_t^{SU,out}\right) \\ &= c^{ASR,in} F^{ASR,in}(L_t) T^{ASR,in} + Q_t^{MS,heat} + Q_t^{SU,heat} \end{aligned} \tag{43}$$

with $Q_t^{MS,heat}$ and $Q_t^{SU,heat}$ chosen as decision variables, eliminating bilinear terms such as $F_t^{MS} T_t^{MS,out}$ and $F_t^{SU} T_t^{SU,out}$. Considering that the RIG flow rate is large enough, by ensuring a gap between RIG and MS temperature, such reformulation will not lead to overestimation of heat exchange ability.

For MS-TES heating, constraints (27)-(29) are replaced by

$$0 \leq Q_t^{MS,heat} \leq M b_t^{MS,on}, \tag{44}$$

$$Q_t^{MS,heat} \leq \overline{Q}^{MS,heat}, \tag{45}$$

where $\overline{Q}^{MS,heat}$ is a calibrated upper bound.

Similarly, SUH heating (31)-(33) is represented by

$$0 \leq Q_t^{SU,heat} \leq M\left(b_t^{AS,on} + b_t^{AS,by} + b_t^{AS,cs}\right), \tag{46}$$

$$0 \leq Q_t^{SU,heat} \leq \overline{Q}^{SU,heat}, \tag{47}$$

where $\overline{Q}^{SU,heat}$ is the upper limit on the SUH heating duty, selected similarly to the MS-TES.

Substituting (36), (37), and (41) into (5) yields a linear expression for ROG heat duty:

$$Q_t^{\text{ASR,out}} = c^{\text{ASR,out}} F_t^{\text{ASR,out}} T^{\text{ASR,out}}$$
$$= c^{\text{ASR,out}} \left( \tilde{T}^{\text{ASR,out}} \left( b_t^{\text{AS,on}} F^{\text{ASR,in,0}} + L_t F^{\text{ASR,in,1}} \right) + (b_t^{\text{AS,by}} + b_t^{\text{AS,cs}}) F^{\text{ASR,in,by}} \tilde{T}^{\text{ASR,in,by}} \right), \tag{48}$$

where $F^{\text{ASR,in,0}}$ and $F^{\text{ASR,in,1}}$ are the coefficients obtained from design data.

### 3.4. Hydrogen production and storage

The model for converting renewable energy into electricity and hydrogen is straightforward and has been extensively studied, including in our previous work [5, 46], Therefore, only a brief description is provided here.

**a) Hydrogen production via water electrolysis.** The hydrogen plant in ReP2H systems is typically large in scale. For a project with an annual ammonia output of 200 kt, a 200 MW hydrogen plant is commonly deployed [12], consisting of 40×5 MW electrolyzers. By on-off scheduling of electrolyzers, the plant load can be modeled as continuously adjustable from 0 to 100% [5], and the operating model is formulated as:

$$P_t^{\text{HP}} = c^{\text{HP}} F_t^{\text{HP,H}_2}, \tag{49}$$

$$0 \le F_t^{\text{HP,H}_2} \le \overline{F}^{\text{HP,H}_2}, \tag{50}$$

where $P_t^{\text{HP}}$ and $F^{\text{HP,H}_2}$ are hydrogen production power and flow rate; $c^{\text{HP}}$ is the power-to-hydrogen conversion factor; $\overline{F}^{\text{HP,H}_2}$ is the limit of $F^{\text{HP,H}_2}$.

**b) Hydrogen storage (HS).** Hydrogen storage dynamics are given by mass conservation [46]:

$$n_{t+1}^{\text{H}_2,\text{Buffer}} = n_t^{\text{H}_2,\text{Buffer}} + \left( F_{t+1}^{\text{HP,H}_2} - F_{t+1}^{\text{H}_2,\text{out}} \right) \Delta t, \tag{51}$$

where $n^{\text{H}_2,\text{Buffer}}$ is the hydrogen inventory; $F^{\text{H}_2,\text{out}}$ is the outlet flow, which equals AS inlet hydrogen flow:

$$F^{\text{H}_2,\text{out}} = F^{\text{ASR,H}_2}, \tag{52}$$

Inventory limits and cyclic operation are enforced by:

$$\underline{n}^{\text{H}_2,\text{Buffer}} \le n_t^{\text{H}_2,\text{Buffer}} \le \overline{n}^{\text{H}_2,\text{Buffer}}, \tag{53}$$

$$n_0^{\text{H}_2,\text{Buffer}} = n_\tau^{\text{H}_2,\text{Buffer}}. \tag{54}$$

### 3.5. Battery energy storage (BES)

BES charging and discharging power, $P_t^{\text{BAT,cha}}$ and $P_t^{\text{BAT,dis}}$, are bounded by

$$0 \le P_t^{\text{BAT,dis}} \le b_t^{\text{BAT,dis}} \overline{P}^{\text{BAT,dis}}, \tag{55}$$

$$0 \le P_t^{\text{BAT,cha}} \le b_t^{\text{BAT,cha}} \overline{P}^{\text{BAT,cha}}, \tag{56}$$

where $\overline{P}^{\text{BAT,cha}}$ and $\overline{P}^{\text{BAT,dis}}$ are the limits; $b_t^{\text{BAT,cha}}$ and $b_t^{\text{BAT,dis}}$ are binaries indicating charging and discharging states, following

$$b_t^{\text{BAT,cha}} + b_t^{\text{BAT,dis}} \le 1, \tag{57}$$

$$0 \le b_t^{\text{BAT,cha}}, b_t^{\text{BAT,dis}} \le 1. \tag{58}$$

The state of charge (SOC) must satisfy:

$$E_{t+1}^{\text{BAT}} = E_t^{\text{BAT}} + \left( P_{t+1}^{\text{BAT,cha}} - P_{t+1}^{\text{BAT,dis}} \right) \Delta t, \tag{59}$$

$$\underline{E}^{\text{BAT}} \le E_t^{\text{BAT}} \le \overline{E}^{\text{BAT}}, \tag{60}$$

$$E_0^{\text{BAT}} = E_\tau^{\text{BAT}}. \tag{61}$$

where $E_t^{\text{BAT}}$, $\underline{E}^{\text{BAT}}$, and $\overline{E}^{\text{BAT}}$ are the stored energy and its limits. The typical allowable SOC range is 10% to 90%.

### 3.6. System power balance

In recent commissioned ReP2A projects, renewable power supplies most of the demand, while the grid acts as a backup [2, 28, 29]. In this study, we assume the system can purchase electricity but does not export power, consistent with existing projects [28, 29, 46-51]. Power balance satisfies

$$P_t^{\text{BAT,cha}} - P_t^{\text{BAT,dis}} + P_t^{\text{MS,heat}} + P_t^{\text{SU,heat}} + P_t^{\text{ASR,aux}} + P_t^{\text{HP}} \le P_t^{\text{RES}} + P_t^{\text{Grid}}, \tag{62}$$

$$0 \leq P_t^{\text{Grid}} \leq \overline{P}^{\text{Grid}}. \tag{63}$$

where $P_t^{\text{RES}}$ is the renewable power; $\overline{P}^{\text{Grid}}$ is the grid power purchase limit. Note that the AS auxiliary load and heating power are also included.

### 3.7. Deterministic scheduling model

The scheduling program aims to maximize the overall profit, defined as ammonia sales minus electricity purchase, start-up/shutdown costs, and operation and maintenance (O&M) costs:

$$C^{\text{Profits}} = \sum_{t=1}^{\tau}(-c^{\text{NH}_3}F_t^{\text{ASR,NH}_3} + c^{\text{Grid}}P_t^{\text{Grid}})\Delta t + c^{\text{Startup}}\sum_{t=1}^{\tau}b_t^{\text{AS,Inoff}} - C^{\text{Z}}, \tag{64}$$

where $c^{\text{NH}_3}$ and $c^{\text{Grid}}$ are ammonia and grid electricity prices; $c^{\text{Startup}}$ is the AS start-up cost, which mainly accounts for fatigue and degradation; $C^{\text{Z}}$ is the sum of investment and O&M costs in the scheduling horizon, given by:

$$\begin{cases} C^{\text{Z}} = D(C^{\text{Inv}} + C^{\text{Om}})/365 \\ C^{\text{Inv}} = \sum_{i \in \Omega}\left[k_i E_i^{\max} \frac{\eta(1+\eta)^{L_i}}{(1+\eta)^{L_i} - 1}\right], \\ C^{\text{Om}} = \sum_{i=1}^{N} k_i E_i^{\max} \omega_i \end{cases} \tag{65}$$

where $D$ is the number of days in the scheduling horizon; $C^{\text{Inv}}$ and $C^{\text{Om}}$ denote investment and O&M costs, respectively; $\Omega$ is the set of components; $E_i^{\max}$, $k_i$, and $L_i$ are component rated capacity, per unit investment cost, and service life, taken as 12 years for the BES and 25 years for the others; and $\omega_i$ is the ratio of O&M cost.

To limit ASR temperature fluctuations, a regulation term is also included:

$$C^{\text{Temp}} = \sum_{t=1}^{\tau}\left[(T_t^{\text{ASR}} - T^{\text{ASR,opt}})^2\right], \tag{66}$$

where $T^{\text{ASR,opt}}$ is the design setpoint.

The overall objective is therefore given by

$$J = -w_1 C^{\text{Profits}} + w_2 C^{\text{Temp}}, \tag{67}$$

where $w_1$ and $w_2$ are weights. Economic performance is emphasized by selecting a large $w_1$.

In summary, the deterministic scheduling model is formulated as:

$$\begin{aligned} &\min \ (67), \\ &\text{s.t. } (1)-(2), (4)-(66). \end{aligned} \tag{68}$$

and its solution serves as the baseline for uncertainty analysis.

### 3.8. IGDT-based scheduling under renewable uncertainty

The scheduling model in Section 3.7 is based on a deterministic renewable profile. To account for uncertainty that may cause deviated revenues from expectations, we extend the baseline model using IGDT [52].

Renewable uncertainty is first defined by [40]:

$$U(\alpha, \tilde{P}_t^{\text{WT/PV}}) = \{P_t^{\text{WT/PV}} : \left|\frac{P_t^{\text{WT/PV}} - \tilde{P}_t^{\text{WT/PV}}}{\tilde{P}_t^{\text{WT/PV}}}\right| \leq \alpha\}, \tag{69}$$

where $\tilde{P}_t^{\text{WT/PV}}$ is the baseline wind and solar power; $\alpha$ is the uncertainty level, which represents the maximum deviation from the baseline.

A robust model maximizes the allowable uncertainty level $\alpha_r$ under a specified revenue loss factor $\beta_r$:

$$\begin{aligned} &\max \alpha_r \\ &\text{s.t. } C \geq (1-\beta_r)C_c^{\text{Profits}} \\ &\quad P_t^{\text{WT/PV}} = (1-\alpha_r)\tilde{P}_t^{\text{WT/PV}} \\ &\text{s.t. } (1)-(2), (4)-(64) \end{aligned} \tag{70}$$

where $C_c^{\text{Profits}}$ is the baseline net revenue obtained in (68).

Conversely, an opportunistic model interprets uncertainty as a potential benefit. It minimizes the minimum uncertainty level $\alpha_o$ needed to reach a desired revenue increase by $\beta_o$:

$$\begin{aligned}
&\min \alpha_o \\
&\text{s.t. } C \geq (1+\beta_o)C_c^{\text{Profits}} \\
&\quad P_t^{\text{WT/PV}} = (1+\alpha_o)\tilde{P}_t^{\text{WT/PV}} \\
&\quad \text{s.t. } (1)-(2), (4)-(64)
\end{aligned} \tag{71}$$

The IGDT-based evaluation procedure is as follows. The deterministic model (68) provides baseline revenue $C_c^{\text{Profits}}$; robust and opportunistic models (70) and (71) then evaluate uncertainty margins $\alpha_r$ and $\alpha_o$. A larger $\alpha_r$ in the robust case indicates stronger tolerance to adverse conditions, while a smaller $\alpha_o$ reflects better ability to exploit favorable conditions. Both models remain MILP and are solvable with standard solvers.

The robust solution is used for day-ahead or intraday scheduling to ensure feasibility under worst-case deviations, while the opportunistic solution informs rolling re-optimization when favorable conditions persist.

## 4. Case Studies

To evaluate the effect of integrating MS-TES into the ReP2A system and the proposed scheduling method, we conduct case studies based on a real-world project in northern China. Simulations are implemented in *Wolfram Mathematica 13.0*, and the optimization problem is solved using *Gurobi 12.0.3*.

### 4.1. Case setups

The system includes a 450 MW wind farm and a 150 MW photovoltaic (PV) plant. Historical data shown in Fig. 4 are used as the baseline scenario, and renewable uncertainty is further analyzed using IGDT. The AS section has a rated capacity of 24.9 t/h (≈200 kt/yr), with its operating load between 30% and 110% of the rating. HS capacity is 150,000 Nm$^3$. Two BES configurations are considered in this study, i.e., a large BES (32 MWh/8 MW) and a small BES (4 MWh/1 MW). The grid electricity price is 1.2 CNY/kWh.

Because the studied system is renewable-dominated with low grid dependence, the ammonia sales price is set to 4,200 CNY/t, reflecting a green premium, e.g., maritime fuel, rather than conventional fossil-based ammonia. The ASR startup cost $c^{\text{Startup}}$ is set to 100,000 CNY. With thermal parameters derived in Appendix B, other technical and economic parameters are listed in Tables 1 and 2 in Appendix C.

The scheduling horizon is 15 days with a 1 h time step. This captures both multi-day renewable variability and the thermal inertia of the ASR and MS-TES, enabling a clear assessment of thermal storage effects.

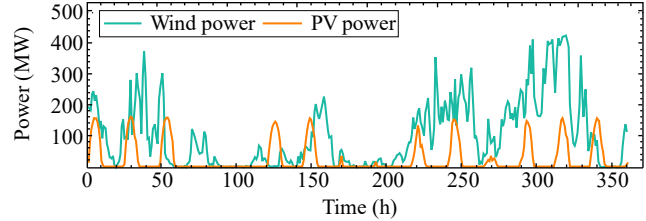

**Fig. 4**. Wind and PV power output profiles in the 15-day baseline scenario.

### 4.2. Comparison benchmarks

Five configurations are compared to assess the value of MS-TES. All share identical inputs and constraints. The configurations of energy, hydrogen, and heat storage modules are summarized in Table 2.

**Scheme 1 (large BES+HS, no MS-TES):** includes HS and a 32 MWh BES. ASR startup and standby heating is provided solely by the SUH. It represents a conventional ReP2A system [25, 39, 43] and serves as the benchmark.

**Scheme 2 (HS only, no BES or MS-TES):** includes HS but no BES or MS-TES. It represents operation without electrical and thermal buffering and is used to examine limitations under prolonged renewable shortages.

**Scheme 3 (MS-TES only, no BES or HS):** includes MS-TES without BES or HS. It isolates the contribution of MS-TES and thus quantifies its standalone impact.

**Scheme 4 (large BES+HS+MS-TES):** adds MS-TES to Scheme 1. It evaluates the effect of simply stacking electrical and thermal storage.

**Scheme 5 (small BES+HS+MS-TES):** replaces the large BES in Scheme 4 with a 4 MWh unit to test coordinated electro-thermal storage at lower cost.

Table 2. Energy, hydrogen, and heat storage options in each comparison scheme

| Component | HS | BES (32MWh) | BES (4MWh) | MS-TES |
|---|---|---|---|---|
| **Scheme 1** | √ | √ | × | × |
| **Scheme 2** | √ | × | × | × |
| **Scheme 3** | × | × | × | √ |
| **Scheme 4** | √ | √ | × | √ |
| **Scheme 5** | √ | × | √ | √ |

Note: √ Included; × Not included.

## *4.3. Scheduling results and comparison*

### *4.3.1. Scheme 1 (32 MWh BES+HS, no MS-TES)*

Scheme 1 serves as the benchmark. Fig. 5 shows the optimal scheduling results. As shown in Figs. 5(a)-(b), hydrogen production closely follows renewable output. BES and HS partially smooth fluctuations but provide a limited duration. HS can support the AS for just over 10 h, as observed for $t$=70-81 h and $t$=259-277 h in Fig. 5(b). BES provides even shorter support, shown in Fig. 5(a) ($t$=39-45 h and $t$=300-305 h).

When low renewable output persists beyond the combined support of BES and HS ($t$=82-110 h in Fig. 5(a)), AS enters hot standby, shown in Fig. 5(c). At this stage, maintaining reactor temperature becomes critical. As shown in Fig. 5(d), the SUH operates frequently once the ASR temperature approaches its lower limit (693 K). Fig. 5(a) shows that the BES is quickly depleted and cannot sustain long-duration heating.

Overall, Scheme 1 can mitigate short-term variability but performs poorly under prolonged energy shortages. Maintaining temperature with BES alone leads to rapid depletion and high cost.

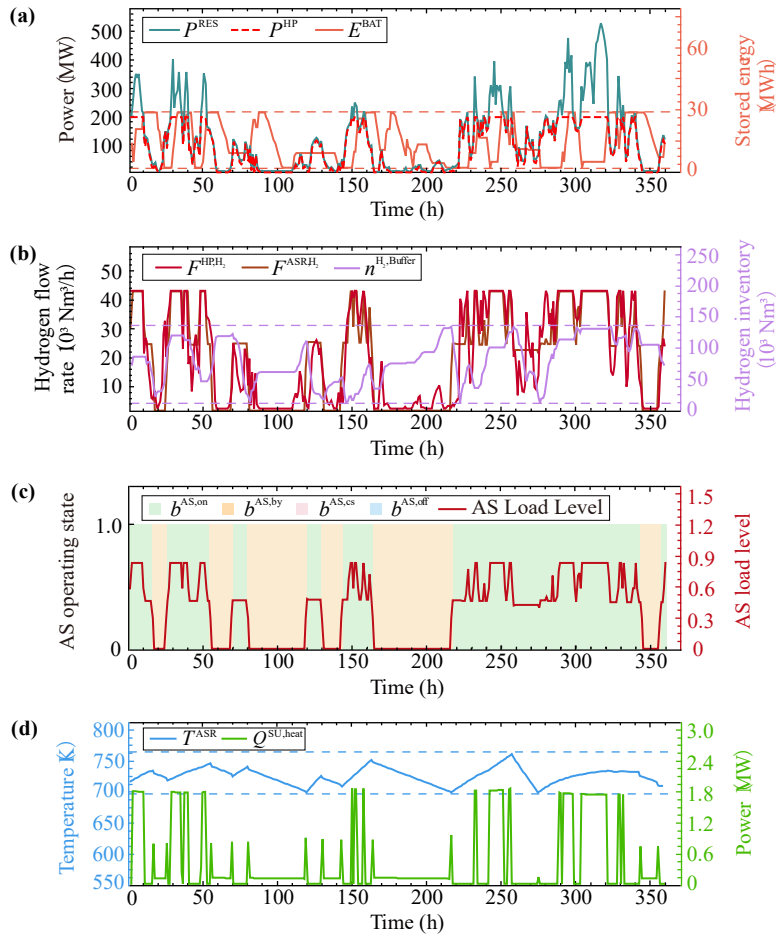

**Fig. 5**. Optimal scheduling results for Scheme 1 (32 MWh BES+HS, no MS-TES). (a) Renewable power generation, hydrogen production power, and the state of BES. (b) Hydrogen production, AS hydrogen consumption, and the state of HS. (c) Operating state and load level of AS. (d) ASR temperature and SUH heat duty.

*4.3.2. Scheme 2 (HS only)*

Scheme 2 excludes both electrical and thermal buffering from BES and MS-TES. Fig. 6 shows that hydrogen production therefore follows renewable output more closely and fluctuates more strongly than Scheme 1.

During extended low-renewable periods (e.g., $t$=163-201 h), the AS cannot remain in hot standby and enters a 38-hour shutdown, shown in Fig. 6(c). The ASR temperature then declines continuously, depicted in Fig. 6(d). When renewable output recovers at around $t$=220 h, production cannot resume immediately; a 12 h cold start is required, during which available power is used for reheating, delaying ammonia production, and the resulting temperature swings intensify reactor thermal cycling. In fact, it leads to the highest start-up and shutdown counts and highest temperature variations among all schemes.

Scheme 2 highlights the limitations of HS alone. Although it provides short-term buffering, it cannot prevent shutdown or large temperature swings during prolonged shortages.

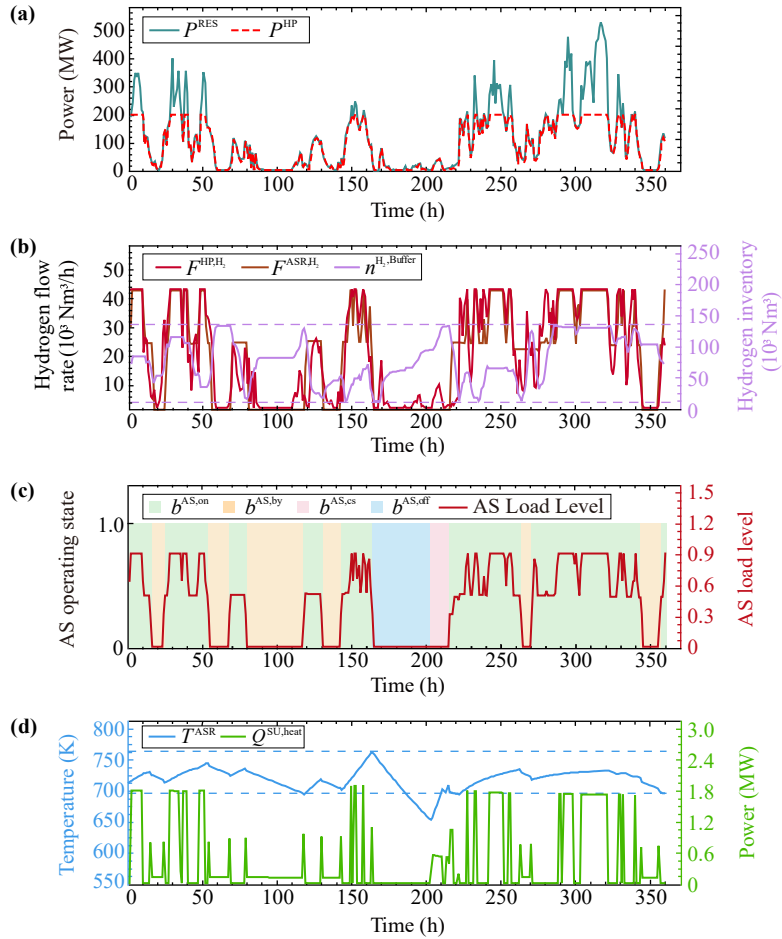

**Fig. 6**. Optimal scheduling results for Scheme 2 (HS only). (a) Renewable power generation and hydrogen production. (b) Hydrogen production, AS hydrogen consumption, and the state of HS. (c) Operating state and load level of AS. (d) ASR temperature and heat duty of SUH.

*4.3.3. Scheme 3 (MS-TES only)*

Scheme 3 isolates the effect of MS-TES. Without BES or HS, hydrogen production directly follows renewable output, as shown in Figs. 7(a), and fluctuations propagate to the AS, shown in Fig. 7(b). Despite this, the AS avoids shutdown and switches only between production and hot standby; see Fig. 7(b). During low-renewable periods (e.g., $t$=43-145 h and $t$=164-231 h), AS remains in hot standby instead of shutting down.

MS-TES stores heat during high-output periods and releases it during energy shortages, as depicted in Fig. 7(c). The SUH plays only a supplementary role. As the result shown in Fig. 7(d), ASR temperature remains smooth across the entire horizon, including prolonged low-renewable intervals. MS-TES also supports production when reaction heat is insufficient to offset thermal losses at $t$=10 h and $t$=239 h.

However, this improved temperature stability comes at the cost of production. As shown in Figs. 7(a) and 7(b), during low-output periods at $t$=43-145 h and $t$=164-231 h, the system remains in hot standby with zero load.

Fig. 8 further illustrates the heat balance under hot standby, at t=270 h. MS-TES supplies heat of 86.2 kW at an outlet temperature of 413.6 °C, which closely matches ASR heat loss of 85.6 kW, maintaining thermal balance and preventing ASR temperature drop.

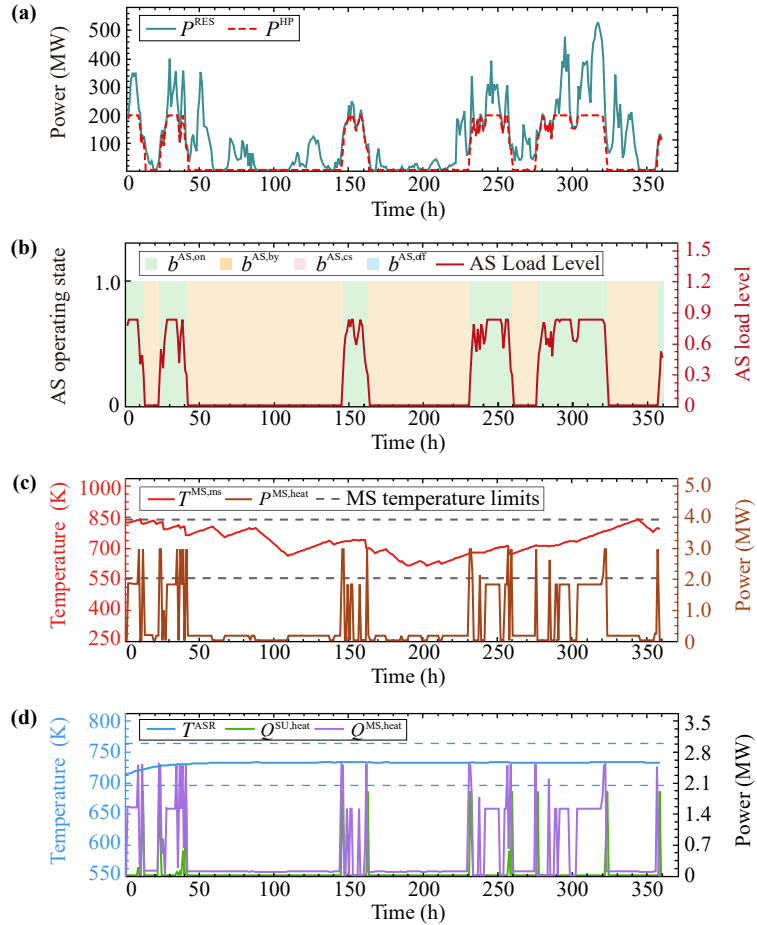

**Fig. 7**. Optimal scheduling results for Scheme 3 (MS-TES only, no BES and HS). (a) Renewable power generation and hydrogen production. (b) AS operating state and load level. (c) MS temperature and electric heating power of the MS-TES. (d) ASR temperature and heat duty of SUH and MS-TES.

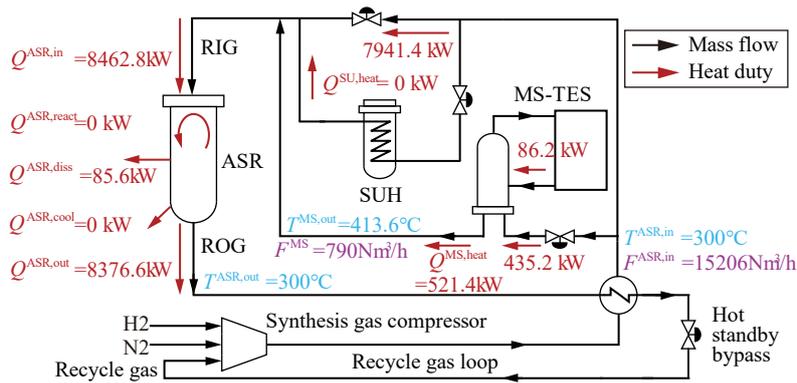

**Fig. 8**. Heat flow in the AS loop under hot standby (at *t*=270 h).

### 4.3.4. Scheme 4 (32 MWh BES+HS+MS-TES)

Scheme 4 combines MS-TES with large BES and HS. Fig. 9(c) shows that the AS operates only in production or hot standby, with no shutdown events. Compared with Scheme 1, ASR temperature fluctuations, in both amplitude and frequency, are significantly reduced. Fig. 9(e) shows that MS-TES provides smoother thermal support, replacing frequent SUH operation.

However, production changes little. The timing of production and standby periods remains similar to Scheme 1, indicating that hydrogen supply remains the limiting factor during prolonged shortages. Meanwhile, the role of the large BES diminishes once MS-TES provides sustained heating.

Overall, Scheme 4 shows that adding MS-TES mainly improves thermal stability but offers limited production gains when large BES and HS are already present.

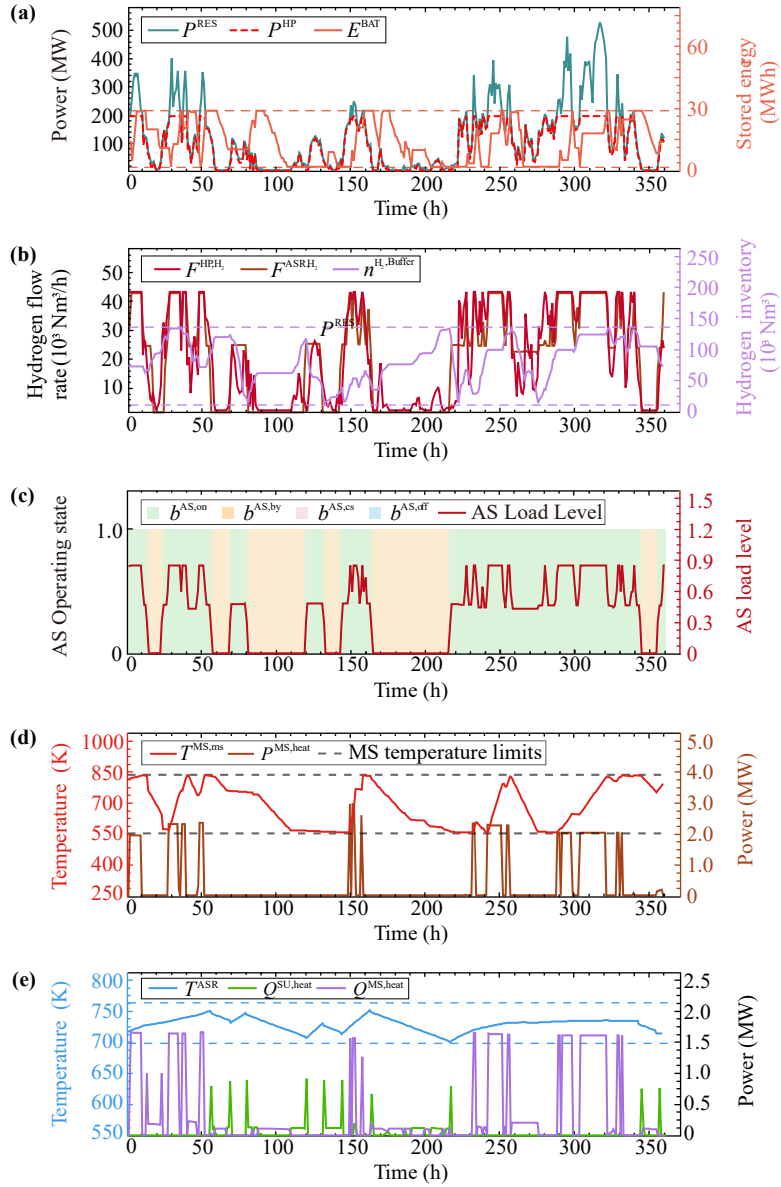

**Fig. 9**. Optimal scheduling results for Scheme 4 (32 MWh BES+HS+MS-TES). (a) Renewable generation, hydrogen production power, and the state of BES. (b) Hydrogen production, AS hydrogen consumption, and the state of HS. (c) AS operating state and load level. (d) MS temperature and electric heating power of the MS-TES. (e) ASR temperature and heat duty of SUH and MS-TES.

*4.3.5. Scheme 5 (4 MWh BES+HS+MS-TES)*

Scheme 5 replaces the large BES in Scheme 4 with a 4 MWh unit, and Fig. 10 shows the scheduling results. Despite the reduced BES capacity, system behavior remains similar to Scheme 4. Hydrogen production, supply, and AS operating states show only minor changes. MS-TES continues to provide most thermal support, with the SUH used only occasionally. ASR temperature remains stable and exhibits smaller fluctuations than in Scheme 1.

Some short low-renewable periods ($t$=271-277 h) lead to temporary hot standby, whereas larger BES configurations maintain production. However, the additional ammonia output from installing a large BES is limited and does not justify its cost.

Overall, Scheme 5 achieves the best balance. Compared with Scheme 2, it avoids shutdown and large temperature swings. Compared with Scheme 1, it further stabilizes ASR temperature. These results indicate that MS-TES can provide long-duration thermal support, allowing a small BES to handle short-term fluctuations at lower cost.

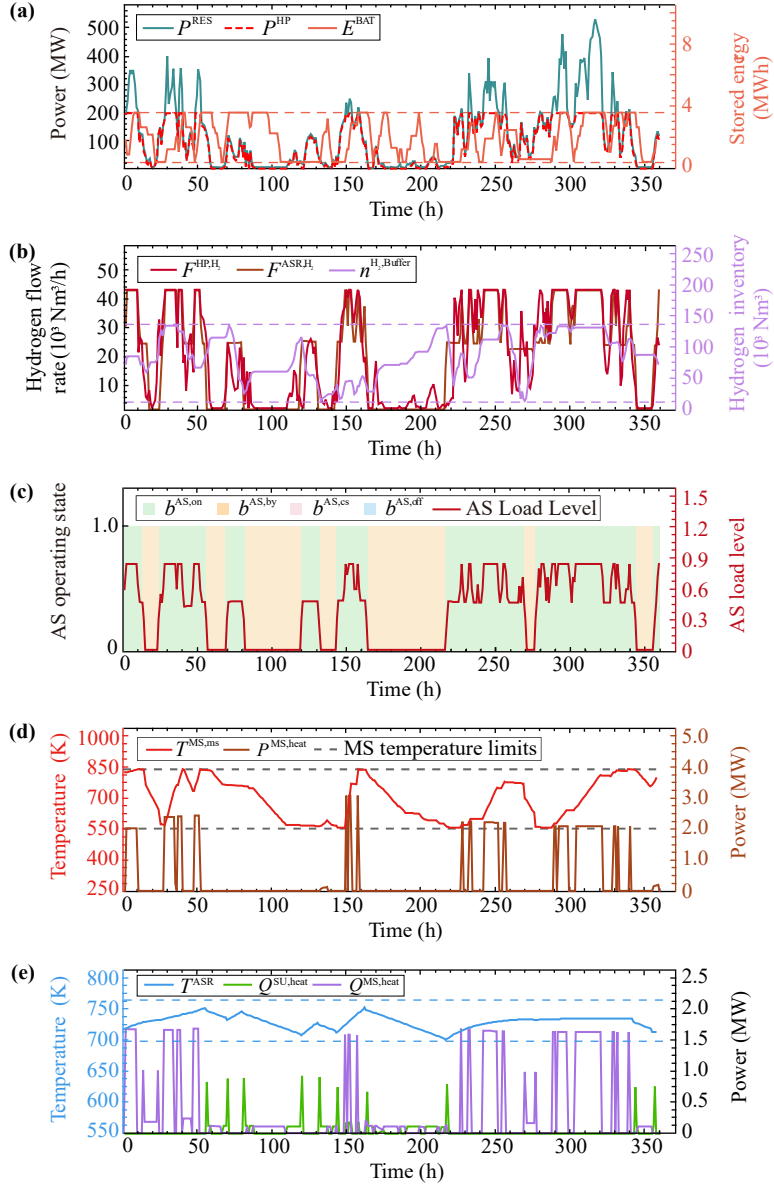

**Fig. 10**. Optimal scheduling results for Scheme 5 (4 MWh BES+HS+MS-TES). (a) Renewable generation, hydrogen production power, and the state of BES. (b) Hydrogen production, AS hydrogen consumption, and the state of HS. (c) AS operating state and load level. (d) MS temperature and electric heating power of the MS-TES. (e) ASR temperature and heat duty of SUH and MS-TES.

*4.3.6. Comprehensive comparison between different schemes*

Table 3 summarizes the scheduling performance of the five schemes and highlights the trade-off between production, thermal stability, and economic performance. Schemes 1-5 are compared in terms of start-up and shutdown counts, ammonia production, costs, net revenue, renewable utilization, and ASR temperature behavior. Because frequent start-up and shutdown increase fatigue and maintenance requirements, they are also included in the cost evaluation. Temperature stability is measured by cumulative temperature variation, defined as the sum of absolute temperature changes between consecutive time steps.

Compared with Scheme 2, Scheme 1 slightly increases ammonia output and renewable utilization, eliminates start-up and shutdown events, and reduces cumulative temperature variation from 526.73 K to 513.55 K. This indicates that the large BES improves operational continuity and temperature stability. However, its higher cost reduces net revenue by 98.3 thousand CNY.

Scheme 3 achieves the most stable temperature profile but produces less ammonia due to the absence of BES and HS, resulting in the highest electricity purchase cost and the lowest net revenue.

Comparing Schemes 1, 4, and 5 clarifies the roles of MS-TES and BES. Relative to Scheme 1, Scheme 4 yields similar ammonia output and net revenue but reduces cumulative temperature variation to 337.75 K, confirming that MS-TES improves temperature stability. Scheme 5 further reduces BES capacity from 32 MWh to 4 MWh. Although ammonia output decreases slightly, it maintains zero start-up and shutdown and similar temperature stability. Its lower cost increases net revenue to 5004.5 thousand CNY, exceeding Schemes 1 and 4 by 101.3 and 96.1 thousand CNY, respectively.

Although Scheme 2 achieves net revenue close to Scheme 5, it shows poorer temperature stability, with start-up and shutdown events and much larger temperature swings. The revenue difference is small because the current objective does not fully monetize thermal fatigue and service life losses. However, if degradation and shortened service life are accounted for, Scheme 2's economic performance is questionable. In contrast, Scheme 5 achieves much better stability with only a 0.43% increase in cost. Overall, Scheme 5 provides the best balance, combining zero start-up and shutdown, low temperature fluctuation, and the highest net revenue in the base 15-day scenario.

Table 3. The scheduling results of the five schemes under the 15-day base-case scenario

| Scheme | AS start-ups and shutdowns | Ammonia yield (t) | Grid electricity purchase cost ($10^4$ CNY) | Investment and O&M costs ($10^4$ CNY) | Net revenue ($10^4$ CNY) | Cumulative temperature variation (K) | Renewable utilization rate (%) |
|---|---|---|---|---|---|---|---|
| Scheme 1 (32 MWh BES+HS, no MS-TES) | 0 | 3441.56 | 3.48 | 951.66 | 490.32 | 513.55 | 81.85 |
| Scheme 2 (HS only) | 2 | 3437.95 | 11.83 | 921.96 | 500.15 | 526.73 | 81.21 |
| Scheme 3 (MS-TES only) | 0 | 2254.59 | 19.05 | 914.36 | 13.51 | 26.66 | 54.80 |
| Scheme 4 (32 MWh BES+HS+MS-TES) | 0 | 3442.24 | 2.99 | 951.91 | 490.84 | 337.75 | 81.93 |
| Scheme 5 (4 MWh BES+HS+MS-TES) | 0 | 3431.70 | 14.94 | 925.93 | 500.45 | 338.10 | 81.41 |

### 4.4. IGDT-based uncertainty evaluation

The IGDT framework in Section 3.8 is used to test whether the preferred storage configuration remains effective when renewable output deviates from the baseline, and the results are shown in Fig. 11. Scheme 3 is excluded because its robust model is infeasible.

As shown in Fig. 11(a), robustness decreases with higher revenue targets for all schemes. Schemes 2 and 5 show the highest robustness, with Scheme 5 consistently performing slightly better. On average, Scheme 5 improves the robustness index by about 2.12% over Scheme 2 and by prominent absolute values of 0.04 and 0.03 over Schemes 1 and 4.

The opportunity curves in Fig. 11(b) show that Schemes 2 and 5 require smaller favorable renewable deviations to achieve a given revenue increase. Scheme 5 reduces the required uncertainty level by 88.53% and 70.56% compared with Schemes 1 and 4. It performs slightly worse than Scheme 2 at low targets but better at high targets, indicating stronger performance under stricter profit requirements.

Overall, the IGDT analysis confirms that Scheme 5 provides both strong robustness and favorable opportunity performance, outperforming conventional large-BES configurations.

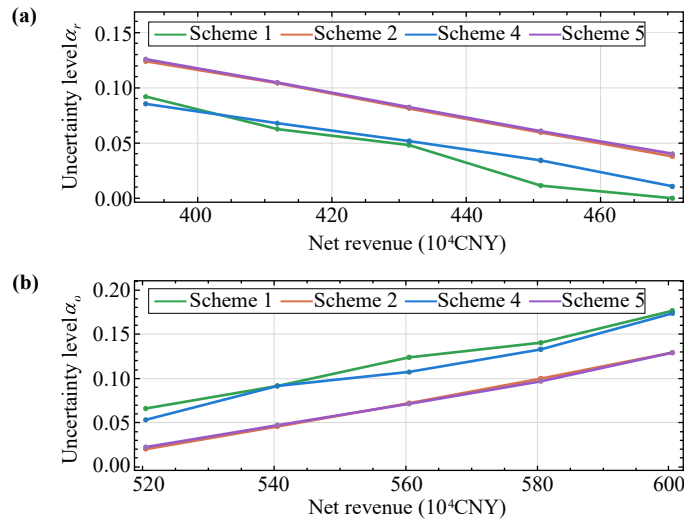

**Fig. 11**. IGDT-based comparison of uncertainty boundaries across different schemes. (a) Robustness. (b) Opportunity.

## 4.5. Sensitivity analysis of electricity, hydrogen, and heat storage sizing

Based on Scheme 5, sensitivity analysis is conducted to examine the effects of BES, HS, and MS-TES sizing. Two-parameter sweeps are performed, i.e., (BES, MS-TES), (MS-TES, HS), and (BES, HS).

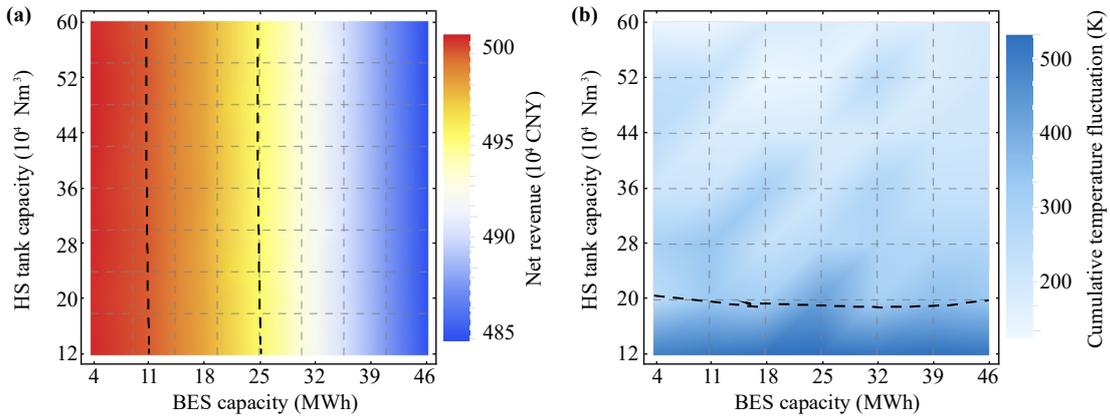

**Fig. 12**. Sensitivity of (a) net revenue, and (b) cumulative temperature variation, to BES capacity and MS-TES sizing.

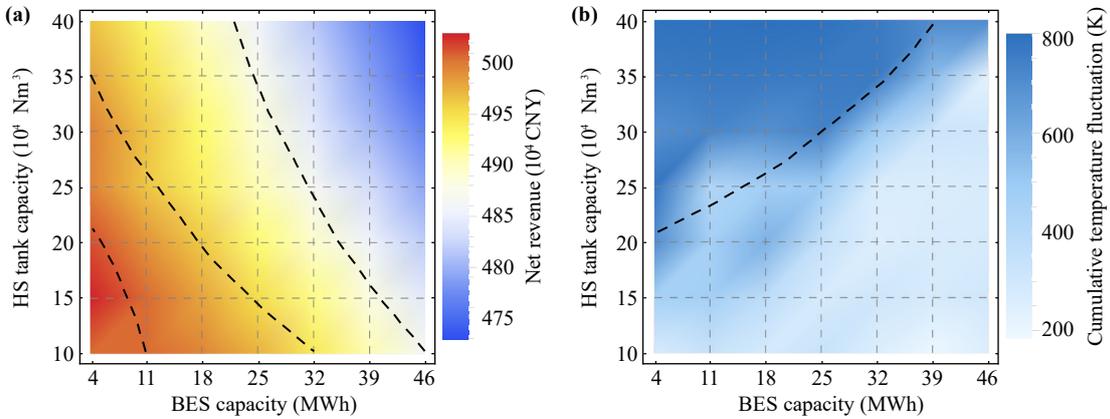

**Fig. 13**. Sensitivity of (a) net revenue, and (b) cumulative temperature variation, to BES and HS capacity.

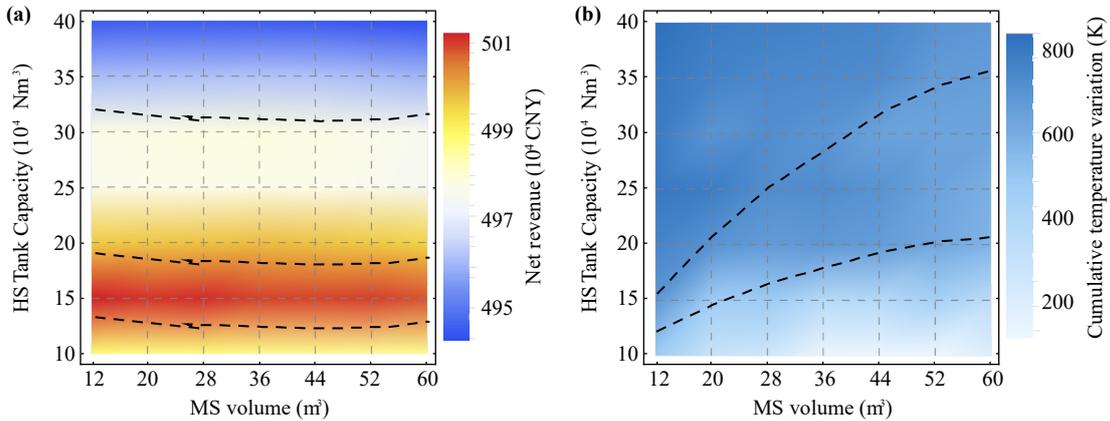

**Fig. 14**. Sensitivity of (a) net revenue, and (b) cumulative temperature variation, to MS-TES sizing and HS capacity.

With HS fixed at 150,000 Nm$^3$, Fig. 12 shows that net revenue decreases mainly with increasing BES capacity, while temperature fluctuation decreases with increasing MS-TES volume. A favorable region emerges with small BES (4-11 MWh) and large MS-TES (36-60 m$^3$), achieving both high revenue and low temperature variation.

With MS-TES fixed at 20 m$^3$, Fig. 13(a) shows net revenue decreases with both BES and HS capacity. The sensitivity to BES capacity is stronger. The net revenue remains above 5.0 million CNY when BES capacity is 4-11 MWh and HS capacity is around 10×10$^4$-15×10$^4$ Nm$^3$, but falls below 4.90 million CNY when BES exceeds 32 MWh and HS exceeds 35×10$^4$ Nm$^3$. Fig. 13(b) shows that temperature variation decreases with BES but increases with HS capacity. In particular, when BES capacity is 4-11 MWh and HS capacity is 30×10$^4$-40×10$^4$ Nm$^3$, cumulative temperature variation exceeds 700 K. This indicates that BES improves thermal stability, whereas excessive HS may worsen thermal stability without sufficient electrical buffering.

With BES fixed, net revenue is more sensitive to HS than to MS-TES. An economic optimum appears at moderate HS capacity (~15×10⁴ Nm³). Temperature fluctuation increases with HS but decreases with MS-TES, showing that MS-TES effectively mitigates thermal variation.

Overall, net revenue is driven mainly by BES and HS sizing, while MS-TES primarily improves temperature stability. Optimal performance requires coordinated sizing rather than expansion of a single resource.

### 4.6. Year-round validation

Annual wind and PV data (including an extended period of zero PV due to scheduled maintenance), shown in Fig. 15, are used to validate performance across 24 scenarios of 15 days each. Results are summarized in Fig. 16.

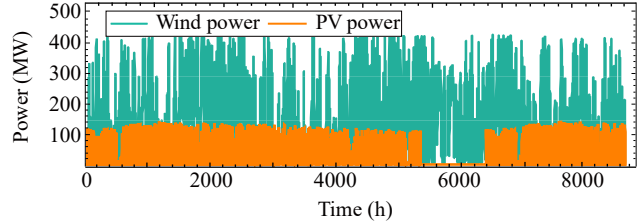

**Fig. 15**. Annual wind and PV generation for year-round validation

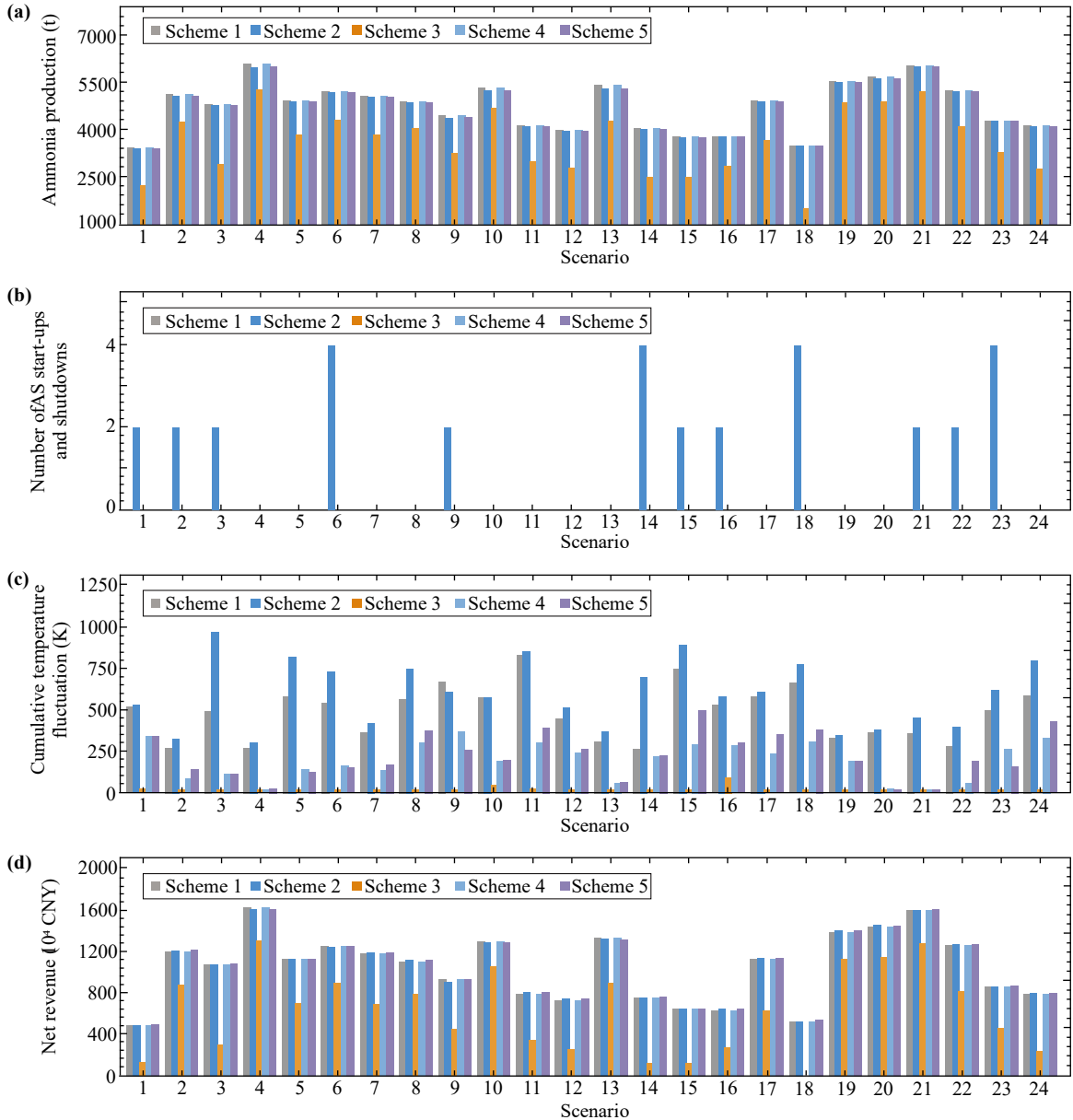

**Fig. 16**. Comparison of scheduling results across a year. (a) Ammonia production. (b) Number of AS start-ups and shutdowns. (c) Cumulative temperature variation. (d) Net revenue.

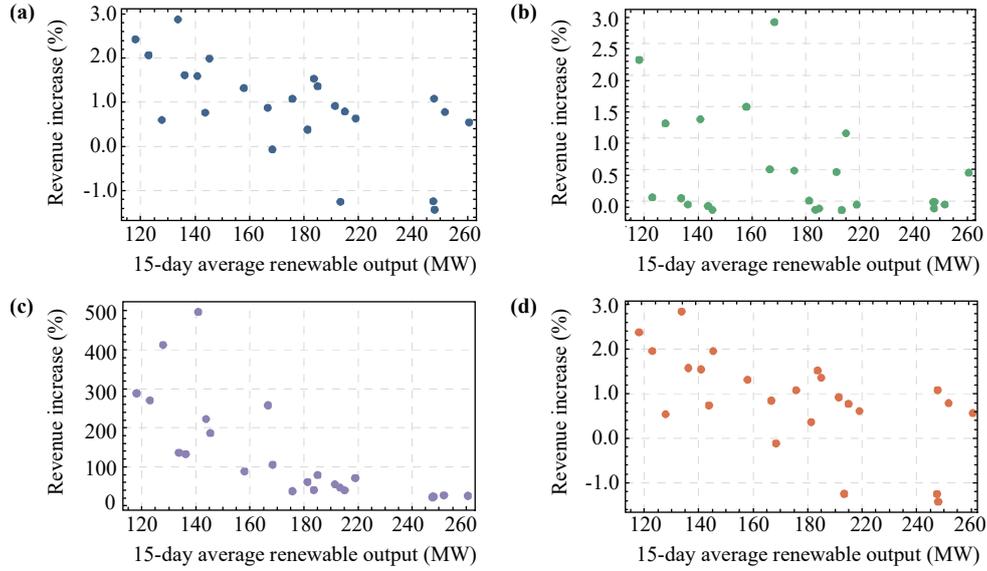

Fig. 17. Net revenue improvement of Scheme 5 relative to (a) Scheme 1, (b) Scheme 2, (c) Scheme 3, and (d) Scheme 4, across 24 scenarios.

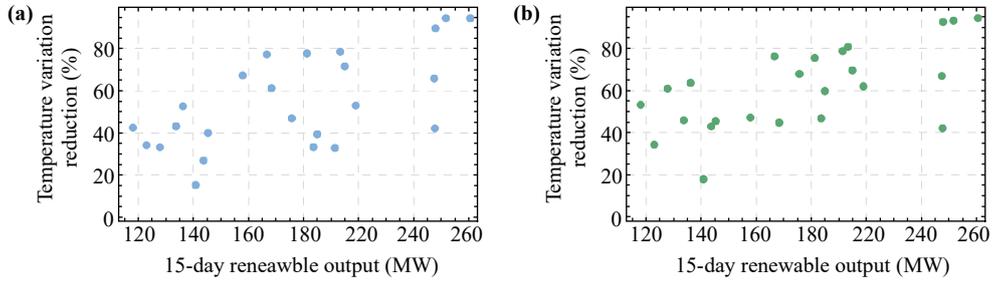

Fig. 18. Cumulative temperature variation reduction of (a) Scheme 5, and (b) Scheme 4, relative to Scheme 1 across 24 scenarios.

Overall, the multi-scenario results are consistent with the base case. Scheme 2 exhibits frequent start-up and shut-down, while the other schemes avoid them. Scheme 3 yields the lowest ammonia output. Schemes 4 and 5 significantly reduce temperature fluctuations compared with Schemes 1 and 2. Scheme 5 achieves the highest net revenue in 20 of 24 scenarios while maintaining low temperature fluctuation and zero start-up and shutdown.

Fig. 17 further shows Scheme 5 (small BES+HS+MS-TES) outperforms Scheme 1 in most scenarios and generally exceeds Scheme 4 (large BES+HS+MS-TES), indicating that the coordinated configuration of a small BES, HS, and MS-TES performs better than large-BES configurations. Gains over Scheme 2 are smaller but positive in more than half the cases. Notably, Scheme 5 outperforms all other configurations in low-renewable cases. However, under resource-rich conditions, its advantage narrows due to limited BES capacity, as the 4 MWh BES cannot fully shift surplus energy.

Fig. 18 confirms that both Schemes 4 and 5 reduce temperature fluctuation relative to Scheme 1, confirming the benefit of MS-TES for ASR temperature stabilization. Scheme 4 is slightly better but at a higher cost.

Table 4. Statistics of annual start-up and shutdown counts and scenario-averaged metrics of all schemes across 24 annual scenarios

| Scheme | Total annual start-up and shutdown counts | Scenario-averaged ammonia production (t) | Scenario-averaged net revenue ($10^4$ CNY) | Scenario-averaged cumulative temperature fluctuation (K) | Scenario-averaged net revenue improvement of Scheme 5 relative to this scheme |
|---|---|---|---|---|---|
| Scheme 1 (32 MWh BES+HS, no MS-TES) | 0 | 4781.11 | 1055.98 | 482.42 | 0.67% |
| Scheme 2 (HS only) | 32 | 4743.34 | 1058.96 | 591.82 | 0.38% |
| Scheme 3 (MS-TES only) | 0 | 3648.45 | 616.31 | 24.72 | 72.48% |
| Scheme 4 (32 MWh BES+HS+MS-TES) | 0 | 4781.94 | 1056.12 | 195.32 | 0.65% |
| Scheme 5 (4 MWh BES+HS+MS-TES) | 0 | 4745.89 | 1063.02 | 224.22 | — |

Table 4 summarizes annual performance. Scheme 5 achieves zero start-up and shutdown, low temperature fluctuation (224.22 K), ammonia output comparable to Schemes 1 and 4, and the highest average net revenue (10.63 million CNY), exceeding Schemes 1-4 by 0.67%, 0.38%, 72.48%, and 0.65%, respectively. The annual results confirm that the coordinated small-BES, HS, and MS-TES configuration is not a case-specific optimum, but a repeatable design choice under diverse renewable conditions.

## 5. Conclusion

This study addresses thermal instability in renewable power-to-ammonia systems by integrating molten-salt thermal energy storage (MS-TES) into the synthesis loop and embedding reactor thermal dynamics into a coordinated scheduling framework. The results highlight three key findings.

1) MS-TES effectively provides long-duration thermal support and directly stabilizes reactor temperature, but it cannot maintain production on its own, as hydrogen availability remains the primary constraint. This confirms that thermal storage complements rather than replaces electrical and material buffering.

2) System performance depends on coordinated resource allocation. A configuration combining small BES, HS, and MS-TES achieves comparable production to large-BES systems while significantly improving thermal stability and reducing overall cost. Once MS-TES is introduced, the marginal benefit of expanding BES becomes limited.

3) The advantages of integrating MS-TES persist under renewable uncertainty. IGDT analysis shows that the coordinated configuration maintains both strong robustness against adverse conditions and improved ability to capture favorable scenarios.

4) Sensitivity analysis further reveals that BES and HS sizing primarily determine economic performance, whereas MS-TES mainly governs thermal stability. Increasing investment in a single resource in isolation is therefore ineffective; coordinated sizing across electricity, hydrogen, and thermal buffering is essential.

Future work will extend the framework to joint capacity planning and operation, incorporate more detailed reactor and catalyst dynamics, and explore real-time control strategies for MS-TES-integrated ammonia synthesis systems.

## Appendix A. Design data of a 200 kt/yr AS process

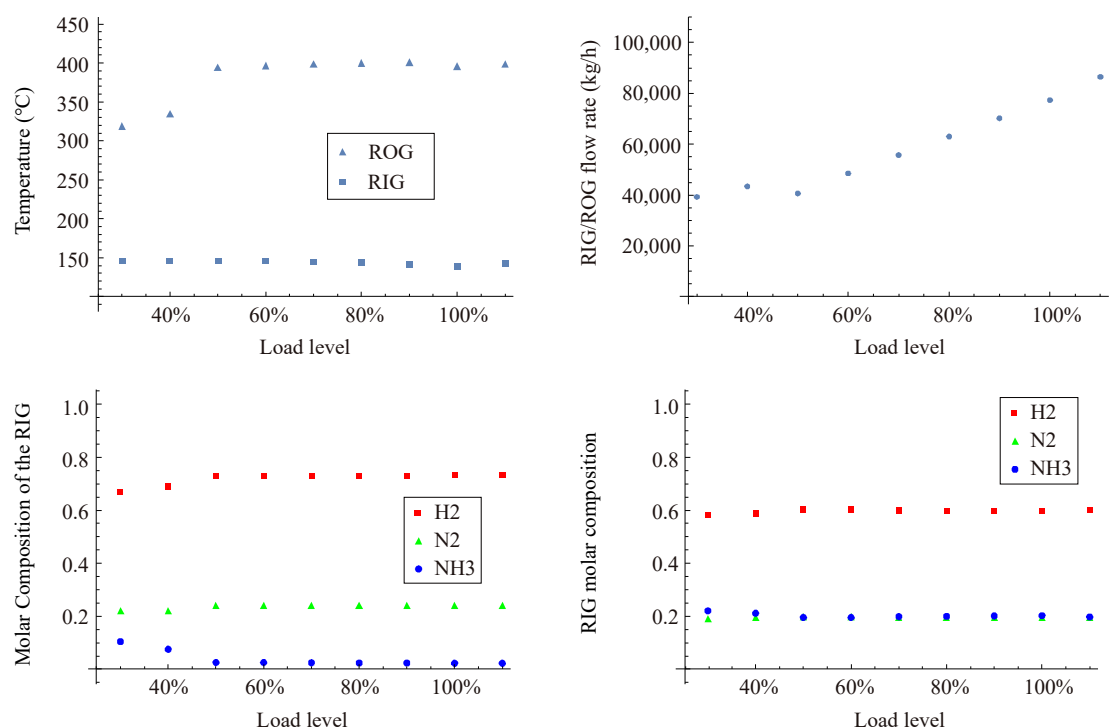

**Fig. A1**. Industrial design data of RIG and ROG flow rate, temperature, and composition for a 200 kt/yr AS process.

**Appendix B. Estimation of thermal model parameters**

To support the simplified thermal model, the ASR thermal capacitance and heat-loss resistance, as well as the heat-loss resistance of the MS-TES, are estimated. The AS unit employs a GC-R301YZ reactor (Nanjing Guochang), with an inner diameter of 1.8 m. The reactor is 20.8 m high, with a catalyst volume of 29.4 m$^3$. The shell is made of Q345R carbon steel (density 7850 kg/m$^3$, specific heat 480 J/(kg·K)). Internal structures are made of 1Cr18Ni9Ti stainless steel (density 7800 kg/m$^3$, specific heat 500 J/(kg·K)) with a total mass of 73 t [55]. The catalyst (A110-1, 87.4 t) is approximated using the specific heat of Fe$_3$O$_4$, its main component.

The total thermal capacitance is the sum of contributions from the shell, internal structures, and catalyst:

$$C^{ASR} = C^{SHELL} + C^{INT} + C^{CAT} = c^{SHELL}m^{SHELL} + c^{INT}m^{INT} + c^{CAT}m^{CAT},$$

where $c^{SHELL/INT/CAT}$ and $m^{SHELL/INT/CAT}$ are the specific heat capacity and mass of the reactor shell/internal structures/catalyst.

The shell volume $V$ is approximated as a hollow cylinder, and the shell mass is estimated as follows:

$$m^{SHELL} = r^{SHELL}V,$$

where $r^{SHELL}$ is the density of the shell material. The corresponding thermal capacitances of the shell, internal structures, and catalyst are 9.849×10$^7$ J/K, 3.65×10$^7$ J/K and 5.681×10$^7$ J/K, respectively, giving the lumped thermal capacitance of 1.918×10$^8$ J/K.

Heat loss from the ASR occurs mainly through conduction across the wall and insulation. These layers are modeled as series thermal resistances, while sidewalls and end surfaces are treated as parallel paths. Radial conduction in the cylindrical structure is approximated using equivalent planar resistances for system-level modeling. The reactor shell (Q345R) has a thermal conductivity of 48 W/(m·K). The insulation layer is high-silica glass fiber (thickness 22 mm, conductivity 0.035 W/(m·K)) [56].

Based on geometry, the outer sidewall area $A^{Side,w}$, single end area $A^{End,w}$, insulation sidewall area $A^{Side,ins}$, and insulation end area $A^{End,ins}$ are 143.759 m$^2$, 3.801 m$^2$, 116.643 m$^2$, and 2.607 m$^2$, respectively. The thermal resistances of the sidewall and one end surface are given by

$$R^{ASR,side/end} = \frac{d^{ASR,w}}{U^{ASR,steel}A^{Side/End,w}} + \frac{d^{ASR,ins}}{U^{ASR,ins}A^{Side/End,ins}}, \tag{B1}$$

where $d^{ASR,w/ins}$ and $U^{ASR,w/ins}$ is the shell/insulation thickness and thermal conductivity. Convective and radiative losses can be implicitly included in the equivalent thermal resistance. The overall equivalent heat-loss resistance of the reactor is then given by:

$$R^{ASR,diss} = \frac{1}{\frac{1}{R^{ASR,side}} + \frac{1}{R^{ASR,end}} + \frac{1}{R^{ASR,end}}} = 0.0052 \text{ K/W}, \tag{B2}$$

The MS-TES uses a binary molten salt (40 wt% KNO$_3$, 60 wt% NaNO$_3$) with density 1924.6 kg/m$^3$ and specific heat capacity 1488 J/(kg·K) [57]. The storage tank has a total volume of 24.8 m$^3$ and contains 20 m$^3$ of molten salt. It is approximately 8.685 m long, with an inner diameter of 2 m and a wall thickness of 10 mm. The tank is made of S32168 stainless steel (thermal conductivity 22.2 W/(m·K)) and insulated with aluminosilicate fiber (thickness 200 mm, conductivity 0.05 W/(m·K) [58]).

Approximating the tank as a cylinder, the outer sidewall area, single end area, insulation sidewall area, and insulation end area are 52.577 m$^2$, 3.205 m$^2$, 66.029 m$^2$, and 4.60 m$^2$, respectively. The thermal resistances of the tank sidewall and one end surface therefore follow

$$R^{MS,side/end} = \frac{d^{MS,w}}{U^{MS,steel}S^{Side/End,w}} + \frac{d^{MS,ins}}{U^{MS,ins}S^{Side/End,ins}}, \tag{B3}$$

where $d^{MS,w}$ and $d^{MS,ins}$ are wall and insulation thickness of the MS tank; $U^{MS,steel}$ and $U^{MS,ins}$ are thermal conductivity of the tank wall and insulation layer. The equivalent heat-loss resistance of the MS-TES module is then calculated as

$$R^{MS,diss} = \frac{1}{\frac{1}{R^{MS,side}} + \frac{1}{R^{MS,end}} + \frac{1}{R^{MS,end}}} = 0.0535 \text{ K/W}. \tag{B4}$$

## Appendix C. Key system parameters for the case study

Table C1. Economic parameters of key components in the ReP2A system

| Equipment | Unit capital cost | Annual O&M cost ratio /% |
|---|---|---|
| Wind power generation equipment | 3700 CNY/kW | 2 |
| Photovoltaic power generation equipment | 3450 CNY/kW | 1 |
| Electrolyzer | 500 CNY/kW | 2 |
| Ammonia synthesis equipment | 1100 CNY/(t/yr) | 2 |
| Hydrogen storage tank | 1750 CNY/kg | 1 |
| Electrochemical battery storage | 1700 CNY/(kWh) | 2 |
| Molten-salt thermal energy storage equipment | 150 CNY/(kWh) | 2 |

Table C2. Key system parameters for the case study

| Parameter | Meaning | Value |
|---|---|---|
| $T^{am}$ | Ambient temperature (K) | 15+273 |
| $c^{ASR,in}$ | RIG specific heat capacity [J/(kg·K)] | 3461 |
| $c^{ASR,out}$ | ROG specific heat capacity [J/(kg·K)] | 3297 |
| $F^{ASR,in,1}$ | RIG flow rate slope coefficient | 13,368 |
| $F^{ASR,in,0}$ | RIG flow rate intercept | 62,661.43 |
| $T^{ASR,in}$ | RIG temperature in production state (K) | 143+273 |
| $T^{ASR,out}$ | ROG temperature in production state (K) | 390.4+273 |
| $\underline{T}^{ASR,act}$ | Lower temperature limit of the ASR in production state (K) | 420+273 |
| $\overline{T}^{ASR,act}$ | Upper temperature limit of the ASR in production state (K) | 490+273 |
| $\overline{Q}^{ASR,cool}$ | Upper limit of the ASR cooling duty (W) | $4.1\times10^6$ |
| $S^{ASR,H_2}$ | Hydrogen consumption at the rated AS load (Nm$^3$/h) | 49,202.4 [46] |
| $S^{ASR,NH_3}$ | Ammonia production rate at rated load (t/h) | 24.9 |
| $F^{ASR,in,by}$ | RIG flow rate in hot standby state (kg/h) | 15,205.89 |
| $c^{HP}$ | Power consumption per unit of hydrogen production (W·h/Nm$^3$) | 4800 [46] |
| $c^{ASR,aux}$ | AS auxiliary power coefficient per unit load rate (W) | $1.59\times10^7$ [46] |
| $P^{ASR,aux,0}$ | Baseline power of AS auxiliary (W) | $3.97\times10^6$ [46] |
| $\eta^{MS,gas,heat}$ | Heat transfer efficiency from MS-TES to the RIG | 0.9 |
| $\eta^{MS,heat}$ | Electric heating efficiency of the MS-TES | 0.95 |
| $r^{ms}$ | MS density (kg/m$^3$) | 1924.6 [53] |
| $c^{ms}$ | MS specific heat capacity [J/(kg·K)] | 1488 [53] |
| $V^{ms}$ | MS volume (m$^3$) | 20 |
| $\overline{P}^{MS,heat}$ | Maximum electric heating power of the MS-TES electric heater (W) | $3.1\times10^6$ |
| $\overline{P}^{SU,heat}$ | Maximum heating power of the SUH (W) | $2\times10^6$ |
| $\underline{n}^{H_2,Buffer}$ | Lower capacity limit of HS tank (Nm$^3$) | $1.5\times10^4$ |
| $\overline{n}^{H_2,Buffer}$ | Upper capacity limit of HS tank (Nm$^3$) | $1.35\times10^5$ |

Note: Parameters without an explicitly cited source in the table are derived from typical design data of a 200,000 t/y ammonia plant.

## Acknowledgement

The authors gratefully acknowledge the financial support from the National Natural Science Foundation of China (52377116 and 52577129).

## Declaration of competing interest

None.

## Data Availability

The data related to this work are available upon request.